\documentclass{article}

\usepackage{microtype}
\usepackage{graphicx}
\usepackage{subfigure}
\usepackage{booktabs} 
\usepackage{multirow}
\usepackage{wrapfig}
\usepackage{threeparttable}

\usepackage{hyperref}
\usepackage{amsthm}
\usepackage{amsmath}
\usepackage{amsfonts}
\usepackage{epstopdf}

\usepackage[accepted]{icml2021}
\def\D{\mathcal{D}}
\def\x{\mathbf{x}}
\def\y{\mathbf{y}}
\def\z{\mathbf{z}}
\def\p{\mathbf{p}}
\def\q{\mathbf{q}}
\def\Z{\mathbf{Z}}
\def\Y{\mathbb{R}^n}
\def\S{\mathcal{S}}
\def\X{\mathcal{X}}
\def\B{\mathcal{B}}
\def\D{\mathcal{D}}
\newtheorem{pro}{Proposition}

\icmltitlerunning{A Value-Function-based Interior-point Method for Non-convex Bi-level Optimization}

\begin{document}
	
	\twocolumn[
	\icmltitle{A Value-Function-based Interior-point Method \\ for Non-convex Bi-level Optimization}

	\begin{icmlauthorlist}
		\icmlauthor{Risheng Liu}{dutru,kl,bz}
		\icmlauthor{Xuan Liu}{dutru,kl}
		\icmlauthor{Xiaoming Yuan}{hk}
		\icmlauthor{Shangzhi Zeng}{hk}
		\icmlauthor{Jin Zhang}{sust,ncam}
	\end{icmlauthorlist}
	\icmlaffiliation{dutru}{DUT-RU International School of Information Science and
		Engineering, Dalian University of Technology.}
	\icmlaffiliation{kl}{Key Laboratory for Ubiquitous Network and Service Software of Liaoning
		Province.}
	\icmlaffiliation{bz}{Pazhou Lab, Guangzhou.}
	\icmlaffiliation{hk}{Department of Mathematics, The University of Hong
		Kong.}
		\icmlaffiliation{sust}{Department of Mathematics, Southern University of Science and Technology.}
		\icmlaffiliation{ncam}{National Center for Applied Mathematics Shenzhen.}
	
	\icmlcorrespondingauthor{Jin Zhang}{zhangj9@sustech.edu.cn}

	\icmlkeywords{Machine Learning, ICML}
	\vskip 0.3in
	]

	\printAffiliationsAndNotice{}  

	\begin{abstract}
		Bi-level optimization model is able to capture a wide range of complex learning tasks with practical interest. Due to the witnessed efficiency in solving bi-level programs, gradient-based methods have gained popularity in the machine learning community.
	In this work, we propose a new gradient-based solution scheme, namely, the Bi-level Value-Function-based Interior-point
	Method (BVFIM). Following the main idea of the log-barrier interior-point scheme, we 
	penalize the regularized value function of the lower level problem into the upper level objective. By further solving a sequence of differentiable unconstrained approximation problems, we consequently derive a sequential programming scheme. 
	The numerical advantage of our scheme relies on the fact that, when gradient methods are applied to solve the approximation problem, we successfully avoid computing any expensive Hessian-vector or Jacobian-vector product. We prove the convergence without requiring any convexity assumption on either the upper level or the lower level objective. Experiments demonstrate the efficiency of the proposed BVFIM on non-convex bi-level problems.
	\end{abstract}
	\section{Introduction}
	In recent years, 
	with the rapid growth of complexity of machine learning tasks, 
	an increasing number of models with hierarchical structures have arisen from various fields~\cite{dempe2018bilevel,liu2021investigating}.
	In these models, according to different hierarchies, parameters are mainly divided into two types,
	for instance, hyper-parameters and model parameters in hyper-parameter optimization~\cite{franceschi2017forward,okuno2018hyperparameter,mackay2018self}, network structures and weights in neural architecture search~\cite{liu2018darts,liang2019darts+,chen2019progressive}, meta learners and base learners in meta learning~\cite{franceschi2018bilevel,rajeswaran2019meta,zugner2018adversarial}, generators and discriminators in adversarial learning~\cite{pfau2016connecting}, actors and critics in reinforcement learning~\cite{yang2019provably}, different components in image processing~\cite{liu2019convergence,ijcai2020-101,liu2020bilevel,liu2020investigating}, etc.
	
	Bi-Level Optimizations (BLOs) have been recognized as important tools to capture these machine learning applications with hierarchical structures.
	Mathematically, BLOs can
	be formulated as the following optimization problem:
	\begin{equation}
		\min\limits_{\x \in \X, \y \in \Y} F(\x,\y),  
		\mathrm{ \ s.t.\  } \y \in {\rm \S}(\x),\label{eq:bilevel}
	\end{equation}
	where ${\rm \S}(\x):=\mathop{\arg\min}_{\y}f(\x,\y)$, the UL constraint $\X \subset \mathbb{R}^m$ is a compact set, and the Upper-Level (UL) objective $F:\X \times\mathbb{R}^n\rightarrow\mathbb{R}$ and the Lower-Level (LL) objective $f:\mathbb{R}^m\times\mathbb{R}^n\rightarrow\mathbb{R}$ are continuously differentiable and jointly continuous functions. Indeed, BLOs are by nature a class of hierarchical problem with optimization problems in the constraints. Specifically, aiming to optimize the UL objective $F$, LL variable $\y$ is selected from the optimal solution set $\S(\x)$ of the LL problem governed by UL variable $\x$. Due to the hierarchical structure and nested relationship between UL and LL problems, BLOs are challenging both theoretically and numerically, especially when the LL solution set $\S(\x)$ is not singleton~\cite{jeroslow1985polynomial}. 

	\subsection{Related Work}
	\begin{table*}[t]
	\caption{Comparing the theoretical results among BVFIM and existing gradient-based methods. $ C^1(\X \times \mathbb{R}^n)$ denotes the set of all continuously differentiable functions on $\X \times \mathbb{R}^n$.} \label{T1}
		\begin{center}
			\begin{small}
				\begin{threeparttable}
					\begin{tabular}{lccccc}
						\toprule
						Category&Methods&LLC&LLS&&Required Conditions\\
						\midrule
						\multirow{6}{*}{EGBM}&\multirow{2}{*}{FHG/RHG}&\multirow{2}{*}{w/}&\multirow{2}{*}{w/}&UL&$F(\x,\y)$ $\in C^1(\X \times \mathbb{R}^n)$.\\
						&&&&LL&$f(\x,\y)$ $\in C^1(\mathbb{R}^m \times \mathbb{R}^n)$.\\
						&\multirow{2}{*}{TRHG}&\multirow{2}{*}{w/}&\multirow{2}{*}{w/}&UL&$F(\x,\y) $ $\in C^1(\X \times \mathbb{R}^n)$ is bounded below.\\
						&&&&LL&$f(\x,\y) $ $\in C^1(\mathbb{R}^m  \times \mathbb{R}^n)$ is $L_f$-smooth and strongly convex.\\
						&\multirow{2}{*}{BDA}&\multirow{2}{*}{w/}&\multirow{2}{*}{w/o}&UL&$F(\x,\y) $ $\in C^1(\X \times \mathbb{R}^n)$ is strongly convex.\\
						&&&&LL&$f(\x,\y) $ $\in C^1(\mathbb{R}^m  \times \mathbb{R}^n)$ is level-bounded in $\y$.\\
						\multirow{2}{*}{IGBM}&\multirow{2}{*}{CG/Neumann}&\multirow{2}{*}{w/}&\multirow{2}{*}{w/}&UL&$F(\x,\y) $ $\in C^1(\X \times \mathbb{R}^n)$.\\
						&&&&LL& $f(\x,\y) $ $\in C^1(\mathbb{R}^m  \times \mathbb{R}^n)$, and $\frac{\partial^2 f(\x,\y)}{\partial \y^2} $ is invertible.\\
						\midrule
						\multirow{2}{*}{IPM}&\multirow{2}{*}{BVFIM}&\multirow{2}{*}{w/o}&\multirow{2}{*}{w/o}&UL&$F(\x,\y) $ $\in C^1(\X \times \mathbb{R}^n)$.\\
						&&&&LL&$f(\x,\y) $ $\in C^1(\mathbb{R}^m  \times \mathbb{R}^n)$.\\
						\bottomrule
					\end{tabular}
				\end{threeparttable}
			\end{small}
		\end{center}
	\end{table*}
	BLOs in Eq.~\eqref{eq:bilevel} are generalizations of several well-known optimization problems as well-noted. Since BLOs can model highly complicated learning tasks as aforementioned, it is not surprising that they are hard to solve.
	In early works, the standard methods for solving BLOs are the random search~\cite{bergstra2012random} which evaluates the learning models through randomly sampled parameter values, or the  more sophisticated Bayesian optimization~\cite{hutter2011sequential}. Unfortunately, these gradient-free methods are only capable of handling more than twenty or so hyper-parameters. Gradient-based
	optimization methods, on the other hand, can handle up to hundreds of hyper-parameters. Specifically, the gradient-based approaches transform the BLO into a Single-Level Optimization (SLO) problem and solve it by using gradient-type methods. To reformulate the BLO into an SLO problem and calculate the gradients, one can replace the LL with a dynamics iteration and derive the gradient by Automatic Differentiation (AD), or apply the implicit function theorem to the optimality conditions of the LL
	problem. According to such different features of calculations, these gradient-based methods
	for BLOs can be classified into two main categories: explicit and implicit. For theoretical support, they both rely on the singleton of LL solution set $\S(\x)$ (a.k.a. LLS), which is restrictive in applications. For numerical implementation, 
	they both require repeated products of vectors and matrices, which are extraordinarily costly in computation power. In response to these two limitations, 
	by reformulating the BLO into an SLO problem via value function approach~\cite{outrata1990on,ye1995optimality}, using the smoothing technique on value function via regularization, and penalizing the smoothed value function to the UL objective via the log-barrier penalty, this paper proposes a value-function-based interior-point gradient-type method. 

	\textbf{Explicit Gradient-Based Methods (EGBMs).}
	The key idea underlying this type of approaches is to hierarchically calculate gradients of UL and LL objectives. 
	Specifically, under the LLS and Lower Level Convexity (LLC), the works in~\cite{franceschi2017forward,franceschi2018bilevel} first calculate gradient representations of the LL objective and then perform either reverse or forward gradient computations (termed as Reverse Hyper-Gradient (RHG) and Forward Hyper-Gradient (FHG)) for the UL sub-problem. In order to address the LLS restriction, under the LLC,~\cite{liu2020generic} proposes Bilevel Descent Aggregation (BDA) which characterizes an aggregation of both the LL and the UL descent information. 
	However, since the EGBMs method requires calculating AD for the entire trajectory of the dynamic iteration of LL, the computation load is heavy to calculate the gradient with reasonable preciseness. In order to reduce the amount of computation,~\cite{shaban2019truncated} proposes Truncated Reverse Hyper-Gradient (TRHG) to truncate the gradient trajectory. However, the efficiency of TRHG is sensitive to the truncated path length. Small truncated path length may deteriorate the accuracy of the calculated gradient and large truncated path length cannot reduce the computation cost. Although \cite{liu2018darts} uses the difference of vectors to approximate gradient, there is no theoretical justification for such type of approximation.

	\textbf{Implicit Gradient-Based Methods (IGBMs).} To reduce the computational burden, another method is to decouple the calculation process of the UL gradient from the dynamic system. For this purpose, IGBMs, which
	people also refer to as implicit differentiation~\cite{pedregosa2016hyperparameter,rajeswaran2019meta,lorraine2020optimizing}, 
	replace the LL sub-problem with an implicit equation. Specifically, when the LL problem is strongly convex (of course in this case the LLS is met), 
	taking advantage of the celebrated implicit function theorem, the gradient of the UL objective
	is calculated by implicit differential equations. However, this scheme needs to repeatedly compute the inverse of the Hessian matrix. Although in practice, the Conjugate Gradient (CG) method or Neumann method are involved for fast inverse computation, however, repeated products of vectors and matrices are still in the core. Therefore, it is still expensive to compute and leads to numerical instabilities, especially when the system is ill-conditioned. 
	Although,~\cite{rajeswaran2019meta} proposes adding a large quadratic term on the LL objective to eliminate the ill-conditionedness, this approach may change the solution set of the BLO problem and cannot generate solutions with desired properties.

	Mathematical requirements for the mentioned methods are listed in Table~\ref{T1}. Indeed, as shown in Table~\ref{T1}, for the theoretical guarantee, except BDA, the LLS and LLC assumptions are required by
	both of the two classes, which are actually too tough to be satisfied in real-world complex tasks.
	Recently,~\cite{Ji2020,Ji2021} give a comprehensive study on the nonasymptotic convergence properties of both implicit and explicit gradient-based methods. Specifically, in~\cite{Ji2021}, the authors provide lower complexity bounds for gradient-type methods and propose accelerated bi-level algorithms to achieve such optimal complexity.
	\subsection{Our Contributions}
	In this paper, to address the issues shared by both the explicit and implicit methods, i.e.,
	restrictive LLS and LLC assumptions in theory and repeated Hessian- and Jacobian-vector products in numerical calculation,
	we initialize a new BLOs solution scheme called Bi-level Value-Function-based Interior-point Method (BVFIM). 
	Our starting point is that the inner Simple Bi-level Optimization (SBO) sub-problem in Eq. \eqref{eq:SBO} can be approximated by more straightforward strategies, i.e., value function approach, a smoothing technique based on regularization and the log-barrier penalization. In particular, we begin with a reformulation of the inner SBO sub-problem by transforming the LL problem solution set constraint into an inequality constraint, using the value function of the regularized LL problem.
	Motivated by an idea behind the Interior-point Penalty Method (IPM),
	we consider a log-barrier function to penalize the value function constraint to the objective of the inner SBO sub-problem. The log-barrier penalty together with regularization further results in a series of SLO problems, thereby approximating the inner SBO sub-problem. Note that thanks to the regularizations, the approximation sub-problems are standard unconstrained differentiable optimization problems, which are therefore numerically trackable by popular gradient-type methods. 
	In the absence of the LLS and LLC assumptions, 
	we prove that BVFIM converges to a true solution of the original BLO.
	Table~\ref{T1} compares the convergence results of BVFIM and the existing methods. It can
	be seen that, for the LL sub-problem, assumptions required in
	previous methods are essentially more restrictive than
	that in BVFIM. More importantly, when solving BLOs without LLS and LLC, no theoretical results can be obtained for these classical methods while BVFIM can still obtain the
	same convergence properties as that in LLS and LLC scenario. 
	Now, we briefly streamline the novelty of our approach and contributions.
	\begin{enumerate} 

		\setlength{\itemsep}{2pt}
		\setlength{\parsep}{0pt}
		\setlength{\parskip}{0pt}
	\item  By penalizing the value function of the LL problem through the log-barrier penalty, together with a regularization smoothing technique, BVFIM constructs a series of unconstrained smooth approximated SLO sub-problems which are numerically trackable by gradient-type optimization methods.

	\item  Compared with the existing gradient-based schemes, the theoretical validity of the proposed BVFIM does not rely on the restrictive LLS and LLC assumptions, 
	thereby capturing a much wider range of bi-level learning tasks.
	\item  
	The proposed BVFIM avoids the computationally expensive Hessian-vector and Jacobian-vector products computation. Therefore, BVFIM can effectively reduce the numerical cost, especially when the LL problem is of large scale.
	 \end{enumerate}

\section{ Gradient-based Methods for BLOs }
	
	In this section, from the optimistic bi-level viewpoint, we offer a unified framework that contains existing implicit and explicit gradient methods as special cases. 
	The framework is also used to analyze our BVFIM later on. Due to the unities of algorithmic forms, the superiority of the proposed BVFIM is clearly observed.

	To address the LLS restriction, we regard BLOs from an optimistic bi-level perspective. Specifically, we transform BLOs in Eq.~\eqref{eq:bilevel} into the following form:
	\begin{equation}
	\min\limits_{\x \in \X}\varphi(\x), \label{eq:op_bilevel}
	\end{equation}
	where $\varphi(\x)$ is the value function of the inner SBO sub-problem, i.e.,
	\begin{equation}
	\varphi(\x):=\inf\limits_{\y \in \S(\x)}F(\x,\y) \quad \mathrm{with} \ \ {\rm \S}(\x):=\mathop{\arg\min}_{\y}f(\x,\y). \label{eq:SBO}
	\end{equation}

	However, $\varphi(\x)$ is a value function of the inner SBO problem,
	by nature it is very ill-conditioned, namely, non-smooth, non-convex and usually with jumps. 
	We next present existing several types of approximations for the gradient (grad.~for short) of $\varphi(\x)$ in a unified manner, i.e., 
	\begin{equation}
		\overbrace{\frac{\partial\varphi(\x)}{\partial \x}}^{\mbox{grad.\ of\ $\x$}} = \overbrace{\frac{\partial F(\x,\y)}{\partial \x} }^{\mbox{direct\ grad.\ of\ $\x$}} + \overbrace{G(\x).}^{\mbox{indirect\ grad.\ of\ $\x$}}\label{eq:inner_step}
	\end{equation}
	Gradient-based methods, either implicit or explicit, access 
	the gradient of $\varphi(\x)$ by approximation with Hessian- and Jacobian-vector products, which results in additional theoretical and computational burdens. 
	Based on the different calculation methods for $G(\x)$, we classify the existing gradient-based methods into EGBMs and IGBMs.

	\textbf{EGBMs.} With an initialization point $\y_0$, a sequence $\{\y_t(\x)\}^{T}_{t=0}$ parameterized by $\x$ is generated as
	\begin{equation}\label{dynamic}
		\y_{t+1}(\x)=\Phi_{t} (\x,\y_t),\ t=0,\cdots,T-1,
	\end{equation}
	where $\Phi_t(\x,\y_t(\x))$ is a smooth mapping
	that represents the operation performed by the $t$-th step of
	an optimization algorithm for solving the LL problem. For example, when the gradient descent method is considered, $\Phi (\x,\y_t(\x))=\y_t(\x)-s\frac{\partial  f}{\partial \y}(\x,\y_{t-1}(\x))$ where $s>0$ denoted the corresponding step size. Approximating $\varphi(\x)$ by $\varphi_T(\x) := F(\x, \y_T(\x))$, an approximation of $\frac{\partial \varphi(\x)}{\partial \x}$ is given by $\frac{\partial \varphi_T(\x) }{\partial \x}= \frac{\partial F(\x,\y)}{\partial \x}  + \frac{\partial F(\x,\y_T)}{\partial \y}\frac{\partial \y_T(\x)}{\partial \x}$. Specifically, in this setting, the indirect gradient $G(\x)$ is specialized as:

	 \begin{equation}\label{requ}
	 G(\x)=\left( \frac{\partial \y_T(\x)}{\partial \x}\right)^\top \frac{\partial F(\x,\y_T)}{\partial \y},
	 \end{equation} 
	 where

	\begin{equation}\small\label{EGBM}
		\frac{\partial \y_t(\x)}{\partial \x}=\left( \frac{\partial \y_{t-1}(\x)}{\partial \x}\right) ^\top\frac{\partial \Phi(\x,\y_{t-1})}{\partial \y}+\frac{\partial \Phi(\x,\y_{t-1})}{\partial \x}.
	\end{equation}

	To ensure the quality of the estimation of $G(\x)$, EGBMs need LLC, most of them need LLS (except BDA), and the number of iterations $T$ should be large enough. Note that Eq.~\eqref{EGBM} contains Hessian- and Jacobian-vector products  when $ \Phi $ denotes the gradient descent operator. 
	Thus computing $G(\x)$ is time-consuming due to Eq.~\eqref{EGBM}.
	
	\textbf{IGBMs.}  
	Under LLS assumption, $\varphi(\x) = F(\x, \y^*(\x))$ with unique $\y^*(\x) \in \mathcal{S}(\x)$. When $\y^*(\x)$ is a differentiable function with respect to $\x$, the gradient $\frac{\partial \varphi(\x)}{\partial \x}$ can be calculated through chain rule as $\frac{\partial \varphi(\x)}{\partial \x} = \frac{\partial F(\x,\y)}{\partial \x}  + \frac{\partial F(\x,\y^*(\x))}{\partial \y}\frac{\partial \y^*(\x)}{\partial \x}$. To obtain $\frac{\partial \y^*(\x)}{\partial \x}$, assuming that $\frac{\partial^2 f}{\partial \y^2}$ is invertible, implicit function theorem is applied on the optimality condition of the LL problem (namely, $0 = \frac{\partial f(\x,\y^*(\x))}{\partial \y} $), i.e.,
	\begin{equation}
	\frac{\partial \y^*(\x)}{\partial \x} = -\left( \frac{\partial^2 f(\x,\y^*(\x))}{\partial \y^2} \right) ^{-1} \frac{\partial^2 f(\x,\y^*(\x))}{\partial \y\partial \x }.
	\end{equation}
	Then $G(\x)$ is specialized as
	\begin{equation}
		\begin{aligned}
			G(\x)= -\left( \frac{\partial^2  f(\x,\y^*(\x))}{ \partial \y \partial \x}\right)^\top\left( \frac{\partial^2 f(\x,\y^*(\x))}{\partial \y^2} \right) ^{-1}&\\
			 \frac{\partial F(\x,\y^*(\x))}{\partial \y}&. \label{eq:implicit_hyper_gradient}
		\end{aligned}
	\end{equation}
	 Note that for the sake of approximation quality, the LL objective $f$ must be strongly convex (or $\frac{\partial^2  f(\x,\y^*(\x))}{\partial \y^2}$ is invertible). In practice, the calculation of $G(\x)$ given in Eq. \eqref{eq:implicit_hyper_gradient} is based on numerical approximations. In particular, 
	CG or Neumann methods which involve Hessian-vector product computations are applied. Thus, the computation of $G(\x)$ is numerically expensive and the accuracy of the approximation heavily depends on the condition number of Hessian or the strong convexity constant of $f$~\cite{grazzi2020iteration}.

	In general, both EGBMs and IGBMs require LLC and LLS  (except BDA) for theoretical convergence. From a computational point of view, the implementation of IGBMs involve repeating productions of Hessian-vector and Jacobian-vector.

\section{Bi-level Value-Function-based Interior-point Method}\label{secBVFIM}
	In this section, we propose a new algorithm that penalizes the regularized value function of the LL problem into the UL objective, and thus approximates the value function of the inner SBO sub-problem $\varphi(\x)$, called Bi-level Value-Function-based Interior-point Method (BVFIM).

	Focusing on $\varphi(\x)$ which is in general non-smooth, to apply gradient-type methods, 
	we will design a series of differentiable functions to approximate $\varphi(\x)$.
	Recall that $\varphi(\x)$ denotes
	the value function of the following parameterized SBO:
	\begin{equation}
		\min\limits_{\y \in \Y } F(\x,\y),   
		\mathrm{ \ s.t.\  } \y\in\arg\min_{\y} f(\x,\y).\label{eq:sp_bilevel}
	\end{equation}
	We can equivalently reformulate Eq.~\eqref{eq:sp_bilevel} into an SLO problem by using the value function of its LL problem, i.e., $f^*(\x)=\min_{\y}f(\x,\y)$. As a consequence, instead of Eq.~\eqref{eq:sp_bilevel}, we study the following parameterized SLO:
	\begin{equation}
		\min\limits_{\y \in \Y } F(\x,\y),   
		\ \mathrm{ \ s.t.\  }\  f(\x,\y)\leq f^*(\x).\label{eq:con_bilevel}
	\end{equation}
	However, the inequality constraint $f(\x,\y)\leq f^*(\x)$ is ill-posed, in the sense that $f^*(\x)$ is non-smooth and the constraint does not satisfy any standard regularity condition. To circumvent such difficulty, inspired by~\cite{borges2020a}, we relax such constraint by replacing $f^*(\x)$ with the value function of the regularized LL problem, i.e., 
	\begin{equation}
		f_{\mu}^*(\x)=\min\limits_{\y \in \Y}f(\x,\y)+\frac{\mu_1}{2}\Vert \y \Vert^{2}+\mu_2,\label{BVFIM2}
	\end{equation}
	where $\mu = (\mu_1, \mu_2)$, and $\mu_1, \mu_2$ are two positive constants. $\frac{\mu_1}{2}\Vert \y \Vert^{2}$ is introduced for guaranteeing the smoothness of $f_{\mu}^*(\x)$ and $\mu_2$ is for ensuring the feasibility of the relaxed inequality constraint $f(\x,\y) < f_{\mu}^*(\x)$. In the following, $f_{\mu}^*$ is shown to be differentiable under mild conditions.
	Now, penalizing the relax inequality constraint $f(\x,\y)\leq f_{\mu}^*(\x)$ to the objective by the log-barrier penalty gives us the following smooth approximation of $\varphi(\x)$:
	\begin{equation}\small
		\varphi_{\mu,\theta,\tau}\left( \x\right) =\min_{\y \in \Y} F(\x,\y) + \frac{\theta}{2}\Vert \y \Vert^{2}-\tau \ln (f_{\mu}^*(\x)-f(\x,\y)), \label{BVFIM}
	\end{equation}
	where $(\mu,\theta,\tau)>0$. The additional regularized term $\frac{\theta}{2}\Vert \y \Vert^{2}$ is for ensuring the smoothness of $\varphi_{\mu,\theta,\tau}\left( \x\right)$. It will be shown in the next section that $\varphi_{\mu,\theta,\tau}\left( \x\right) \rightarrow \varphi(\x)$ as $(\mu,\theta,\tau) \rightarrow 0$ and $\tau \ln\mu_{2} \rightarrow 0$.
	Now, we show the differentiability of $\varphi_{\mu,\theta,\tau}\left( \x\right)$ and give the formula of its gradient.
	\begin{pro}\label{gradient}
		Suppose $ F(\x, \y) $ and $f(\x, \y)$ are continuously differentiable. Given $\x \in \X$ and $\mu, \theta, \tau > 0$, when
		\begin{equation}\label{y-update111}
			\begin{aligned}
				&\y_{\mu,\theta,\tau}^*(\x) \\
				=\; &\underset{\y \in \Y}{\mathrm{argmin}}~F(\x,\y)+\frac{\theta}{2}\Vert \y \Vert^{2}-\tau \ln (f_{\mu}^*(\x)-f(\x,\y)),
			\end{aligned}
		\end{equation}
		and 
		\begin{equation}
			\z^*_{\mu}(\x)= \underset{\y \in \Y}{\mathrm{argmin}}~f(\x,\y)+\frac{\mu_{1}}{2}\Vert \y \Vert^{2}+\mu_{2},
		\end{equation}
		are unique, then $\varphi_{\mu,\theta,\tau}$ is differentiable and
		\begin{equation}\label{g_varphi}
			\frac{\partial \varphi_{\mu,\theta,\tau}}{\partial \x} \left( \x\right)=\frac{\partial  F(\x,\y_{\mu,\theta,\tau}^*(\x))}{\partial \x}+G(\x),
		\end{equation}
		where 
		\begin{equation}\label{def_G}
			G(\x)=\frac{\tau\left(\frac{\partial f(\x,\y_{\mu,\theta,\tau}^*(\x))}{\partial \x}-\frac{\partial f(\x,\z^*_{\mu}(\x))}{\partial \x} \right)}{f_{\mu}^*(\x)-f(\x,\y_{\mu,\theta,\tau}^*(\x))},
		\end{equation}
		and $f_{\mu}^*(\x) = f(\x,\z^*_{\mu}(\x)) + \frac{\mu_{1}}{2}\|\z^*_{\mu}(\x)\|^2 + \mu_{2}$.
	\end{pro}
	The decreasing proximity from Eq.~\eqref{BVFIM} to the original BLOs leans at the core of our convergence theory, as the regularization sequence $\{(\mu_k,\theta_k,\tau_k)\}$ is tending to zero.	
	In particular, the solutions returned by solving approximate sub-problems $\mathrm{min}_{\x \in \X}\varphi_k(\x)$, where $\varphi_k$ denotes $\varphi_{\mu_k,\theta_k,\tau_k}$,  
	converge to the true solution of BLOs, as $\{(\mu_k,\theta_k,\tau_k)\}$ vanishes with $\tau_k \ln\mu_{k,2} \rightarrow 0$; see Theorem 1 in Section 4.	
	 Naturally, to execute this computing process, we need to calculate the gradient $\frac{\partial  \varphi_{k}}{\partial \x}$, thereby solving each approximate sub-problem efficiently.

	We next illustrate the calculation of $\frac{\partial  \varphi_{k}}{\partial \x}$ at $\x_l$ as a guide for implementation. Indeed, for a user-friendly purpose, this procedure is divided into three steps, according to Proposition~\ref{gradient}.  	
	First, we calculate $\z^*_{\mu_k}(\x_l)$ by solving the regularized LL problem Eq.~\eqref{BVFIM2}, which can be easily done by gradient descent as
	\begin{equation}\label{update_z}
		\z_{k,l}^t = \z_{k,l}^{t-1} - s_1 \left( \frac{\partial  f(\x_l,\z_{k,l}^{t-1})}{\partial \y} + \mu_{k,1} \z_{k,l}^{t-1} \right),
	\end{equation}
	where $s_1 >0$ is an appropriately chosen step size.
	To calculate $\y_{\mu_k,\theta_k,\tau_k}^*(\x_l)$ returned by solving Eq.~\eqref{BVFIM},  the gradient descent scheme reads as
	\begin{equation}\label{update_y}\small
		\begin{aligned}
			\y_{k,l}^t &=\y_{k,l}^{t-1} \\
			&- s_2  \Bigg( \frac{\partial  F(\x_l, \y_{k,l}^{t-1})}{\partial \y}	
			+ \theta_k \y_{k,l}^{t-1}
			+ \frac{\tau_k \frac{\partial  f(\x_l,\y_{k,l}^{t-1})}{\partial \y}}{ f_{k,l}^{T_{\z}} - f(\x_l,\y_{k,l}^{t-1})}  \Bigg),
		\end{aligned}
	\end{equation}
	where $s_2 >0$ is an appropriately chosen step size. In particular, $f_{k,l}^{T_{\z}} = f(\x_l, \z_{k,l}^{T_{\z}}) + \frac{\mu_{k,1}}{2}\|\z_{k,l}^{{T_{\z}}}\|^2 + \mu_{k,2}$	, where 
	$\z_{k,l}^{{T_{\z}}}$ represents the output returned from $T_\z$-step gradient descent in Eq.~\eqref{update_z}.
	
	Following Eqs.~\eqref{g_varphi} and~\eqref{def_G}, we eventually obtain $\frac{\partial  \varphi_{k}(\x_l)}{\partial \x}$ by
	\begin{equation}\label{a_g_vphi}
		\frac{\partial  \varphi_{k}(\x_l) }{\partial \x}\approx \frac{\partial F(\x_l,\y_{k,l}^{T_{\y}})}{\partial \x}+G_{k,l},
	\end{equation}
	with
	\begin{equation}
	G_{k,l} = \frac{\tau_k\left(\frac{\partial f(\x_l,\y_{k,l}^{T_\y})}{\partial \x} - \frac{\partial f(\x_l,\z_{k,l}^{T_\z})}{\partial \x} \right)}{f_{k,l}^{T_\z}-f(\x_l, \y_{k,l}^{{T_{\y}}} )},
	\end{equation}
	where $\y_{k,l}^{{T_{\y}}}$ represnts the output returned from $T_\y$-step gradient descent in Eq.~\eqref{update_y}, and $\z_{k,l}^{T_{\z}}$, $f_{k,l}^{T_\z}$ are stored in the last step.

	\begin{algorithm}
		\caption{Our Solution Strategy for Eq.~\eqref{BVFIM}}\label{alg:ours}
		\begin{algorithmic}[1]
			\REQUIRE $(\mu_k,\theta_k,\tau_k)$, step size $\alpha,s_1,s_2>0$
			
			\FOR {$l=0\rightarrow L$}
			\STATE Calculate $\z^{*}_{\mu_k}(\x_l)$ by $T_\z$-step gradient descent in Eq.~\eqref{update_z}.
			\STATE Calculate $\y_{\mu_k,\theta_k,\tau_k}^*(\x_l)$ by $T_\y$-step gradient descent in Eq.~\eqref{update_y}.
			\STATE Calculate an approximation of $\frac{\partial  \varphi_{k}(\x_l)}{\partial \x}$ by Eq.~\eqref{a_g_vphi}, denoted by $g_l$.
			\STATE $\x_{l+1}=\x_{l}-\alpha g_l$.
			\ENDFOR
			
		\end{algorithmic}  
	\end{algorithm}
	The algorithm to solve approximate sub-problems Eq.~\eqref{BVFIM} is outlined in Algorithm~\ref{alg:ours}.
	Then, as $\{(\mu_k,\theta_k,\tau_k)\}$ vanishes, the solution of Algorithm~\ref{alg:ours} converges to the true solution of Eq.~\eqref{eq:sp_bilevel}.

\section{ Theoretical Investigations}\label{secTheo}
	In this section, we will give the convergence analysis the proposed method. Please notice that all the proofs of our theoretical results are
	stated in the Supplemental Material.

	We first recall an equivalent definition of epiconvergence given in~\cite{AlexanderShapiro2011Perturbation}[page 41].
	\newtheorem{defi}{Definition}
	\begin{defi} \label{def_epic}
		$\varphi_k \stackrel{e}{\longrightarrow} \varphi$ iff for all $\x \in \mathbb{R}^m$ the following two conditions hold:
		\begin{enumerate}
			\setlength{\itemsep}{2pt}
			\setlength{\parsep}{0pt}
			\setlength{\parskip}{0pt}
		\item For any sequence $\{\x_k\}$ converging to $\x$,
		\begin{equation}
			\liminf \limits_{k \rightarrow \infty}\varphi_k(\x_k) \geq \varphi(\x).
		\end{equation}
		\item There is a sequence $\{\x_k\}$ converging to $\x$ such that
		\begin{equation}
			\limsup \limits_{k \rightarrow \infty}\varphi_k(\x_k) \leq \varphi(\x).
		\end{equation}

		\end{enumerate}
	\end{defi}
	For a given function $f(\x, \y)$, we state the property that 
	it is level-bounded in $\y$ locally uniformly in
	$\x \in \X$ in the following
	definition.
	\begin{defi}\label{epicon}
		Given a function $f(\x, \y) : \mathbb{R}^m \times \mathbb{R}^n \rightarrow \mathbb{R}$. If
		for a point $\bar{\x}\in \X \subseteq \mathbb{R}^m$, for any $c \in \mathbb{R}$, there exists $\delta > 0$ along
		with a bounded set $\B \in \mathbb{R}^m$, such that
		\begin{equation}
			\{\y \in \mathbb{R}^n |f(\x,\y)\leq c\} \subseteq \B, ~~ \forall \x \in \B_{\delta}(\bar{\x})\cap \X,
		\end{equation}
		then we call $f(\x, \y)$ is level-bounded in $\y$ locally uniformly
		in $\bar{\x} \in \X$ . If the above property holds for each $\bar{\x} \in \X$, we further call $f(\x, \y)$ is level-bounded in $\y$ locally uniformly in
		$\x \in \X$ .

	\end{defi}

	We will give the convergence result of the proposed algorithm under the following standing assumptions:
	\newtheorem{asu}{Assumption}
	\begin{asu}
		We take the following as our blanket assumption
		\begin{enumerate}
		\setlength{\itemsep}{0pt}
		\setlength{\parsep}{0pt}
		\setlength{\parskip}{0pt}
		\item $\S(\x) $ is nonempty for $\x \in \X$.
		\item  Both $F(\x,\y)$ and $f(\x,\y)$ are jointly continuous and continuously differentiable.
		\item  Either $F(\x, \y)$ or $f(\x, \y)$ is level-bounded in $\y$ locally uniformly in $\x \in \X$.
		\end{enumerate} 
	\end{asu}
	
	To prove the convergence result, we show that $\varphi_k(\x)+\delta_{\X}(\x)\stackrel{e}{\longrightarrow}\varphi(\x)+\delta_{\X}(\x)$, where $\delta_\X(\x)$ denotes the indicator function of set $\X$, i.e., $\delta_\X(\x) = 0 $ if $\x \in \X$ and $\delta_\X(\x) = + \infty$ if $\x \notin \X$. To do this, we need to verify two conditions given in Definition~\ref{def_epic}. 
	
	First, to show the condition 1 in Definition~\ref{def_epic},
	we introduce the
	value function $\psi_{\mu}(\x)$ of the relaxed problem of Eq.~\eqref{eq:con_bilevel}:
	\begin{equation}
		\psi_{\mu}(\x)=\min\limits_{\y \in \Y } F(\x,\y),   
		\ \mathrm{ \ s.t.\  }\  -1 \le f(\x,\y)-f_{\mu}^*(\x)\leq 0. \label{low_relax}
	\end{equation}
	And we propose the following two lemmas.
	\newtheorem{lemma}{Lemma}
	\begin{lemma}\label{lem2}
		Let $\{\mu_k\}$ be a positive sequence such that $\mu_k \rightarrow 0$. Then for any sequence $\{\x_k\}$ converging to $\bar{\x}$,
		\[
		\limsup_{k \rightarrow \infty} f_{\mu_k}^*(\x_k) \le f^*(\bar{\x}).
		\]
	\end{lemma}
	
	\begin{lemma}\label{lem3}
		Given $\x \in \X$, suppose either $F(\x, \y)$ or $f(\x, \y)$ is level-bounded in $\y$ locally uniformly in $\x$.
		Let $\{ \mu_k\} $ be a positive sequence such that $\mu_k \rightarrow 0$, and then for any sequence $\{ \x_k\}$ converging to $\x$,
		\begin{equation}
			\liminf\limits_{k\rightarrow\infty}\psi_{\mu_k}(\x_k)\geq \varphi(\x).
		\end{equation}
	\end{lemma}
	Next, we verify condition 2 in Definition~\ref{def_epic}.
	\begin{lemma}\label{lem4}
		Let $\{(\mu_k, \theta_k, \tau_k)\}$ be a positive sequence such that $(\mu_k, \theta_k, \tau_k) \rightarrow 0$ and $\tau_k \ln\mu_{k,2} \rightarrow 0$. Then for any $\x \in \X$,
		\[
		\limsup_{k \rightarrow \infty} \varphi_k(\x) \le \varphi(\x).
		\]
	\end{lemma}
	
	Now, by combining the results obtained above,
	we can obtain the desired epiconvergence result.
	\begin{pro}\label{prop2}
		Suppose either $F(\x, \y)$ or $f(\x, \y)$ is level-bounded in $\y$ locally uniformly in $\x \in \X$. Let $\{(\mu_k, \theta_k, \tau_k)\}$ be a positive sequence such that $( \mu_k, \theta_k, \tau_k) \rightarrow 0$ and $\tau_k \ln\mu_{k,2} \rightarrow 0$, and then
		\begin{equation}
			\varphi_k(\x)+\delta_{\X}(\x)\stackrel{e}{\longrightarrow}\varphi(\x)+\delta_{\X}(\x).
		\end{equation}
	\end{pro}
	In summary, 
	we can get the following convergence result:
	\newtheorem{thm2}{Theorem}
	\begin{thm2}\label{the1}
		Suppose either $F(\x, \y) $ or $f(\x, \y)$ is level-bounded in $\y$ locally uniformly in $\x \in \X$. Let
		$\{(\mu_k, \theta_k, \tau_k)\}$ be a positive sequence such that $( \mu_k, \theta_k, \tau_k) \rightarrow 0 $ and $\tau_k \ln\mu_{k,2} \rightarrow 0$. Then
				\begin{enumerate}
			\setlength{\itemsep}{0pt}
			\setlength{\parsep}{0pt}
			\setlength{\parskip}{0pt}
			\item We have the following inequality
			\begin{equation}
			\limsup \limits_{k \rightarrow \infty}\left(\inf \limits_{\x \in \X}\varphi_k(\x) \right)\leq \inf \limits_{\x \in \X}\varphi(\x).
			\end{equation}
			\item  If $\x_\ell \in \mathrm{argmin}_{\x \in \X}\varphi_\ell(\x)$, for some sequence $\{\ell\} \subset \mathbb{N} $ and $\x_\ell$ converges to $\tilde{\x}$, then $\tilde{\x} \in \mathrm{argmin}_{\x\in \X} \varphi(\x)$
			and
			\begin{equation}
			\lim \limits_{\ell \rightarrow \infty}\left (\inf \limits_{\x \in \X}\varphi_\ell(\x) \right)= \inf \limits_{\x \in \X}\varphi(\x).
			\end{equation}
		\end{enumerate} 
	\end{thm2}
	\begin{table*}[t]
	\caption{Complexity analysis of various bi-level methods on unconstrained problems in calculating $\frac{\partial\varphi(\x)}{\partial \x}$. We show the key update steps. $\p\in \Y$  and $\q\in \mathbb{R}^m$ are intermediate variables and $\Z\in\mathbb{R}^{m\times n}$ is intermediate matrix. $\mathbf{I}$ is the identity matrix. FHG and RHG calculate Eq.~\eqref{requ} by forward or reverse AD. CG calculates Hessian inverse by solving linear equation, while Neumann approximates Hessian inverse by calculating the limit of the sequence. Please see ~\cite{franceschi2017forward,pedregosa2016hyperparameter,lorraine2020optimizing} for more calculation details.  Note that our method avoid calculating the Hessian- and Jacobian-vector products.  } 
	\begin{center}
		\begin{small}
			\begin{tabular}{lccc}
				\toprule
				Method&\multicolumn{1}{c}{$G(\x)$}&Time&Space\\
				\midrule
				\multirow{2}{*}{FHG}&\multirow{2}{*}{$ \Z_{T}^\top \frac{\partial F(\x,\y_T)}{\partial \y}$ with $\Z_{t}= \frac{\partial^2 f}{\partial \y^2}\Z_{t-1} +  \frac{\partial^2 f}{\partial \y\partial \x }$}&\multirow{2}{*}{$O(cmT)$}&\multirow{2}{*}{$O(mn)$}\\
				&&&\\
				\multirow{2}{*}{RHG}&\multirow{2}{*}{$\q_{-1}${ with }$ \q_{t-1}=\q_{t}+ \left( \frac{\partial^2 f}{\partial \x \partial \y}\right)^\top\p_t,\p_{t-1}=\left( \frac{\partial^2 f}{\partial \y^2}\right)^\top  \p_{t}$}&\multirow{2}{*}{$O(cT)$}&\multirow{2}{*}{$O(m+nT)$}\\
				&&&\\
				\multirow{2}{*}{CG}&\multirow{2}{*}{$-\left( \frac{\partial^2  f(\x,\y_T)}{ \partial \y \partial \x}\right)^\top\q $ with $ \frac{\partial^2 f}{\partial \y^2} \q=\frac{\partial F}{\partial \y} $}&\multirow{2}{*}{$O(c(T+J))$}&\multirow{2}{*}{$O(m+n)$}\\
				&&&\\
				\multirow{2}{*}{Neumann}&\multirow{2}{*}{$-\left( \frac{\partial^2  f(\x,\y_T)}{ \partial \y \partial \x}\right)^\top\q\frac{\partial F}{\partial \y}$ with $\q=\sum_{j=0}^{J}\left(\mathbf{I}-\frac{\partial^{2} f}{\partial \mathbf{y}^2}\right)^{j}$}&\multirow{2}{*}{$O(c(T+J))$}&\multirow{2}{*}{$O(m+n)$} \\
				&&&\\
				\midrule
				\multirow{3}{*}{BVFIM}&\multicolumn{1}{c}{\multirow{3}{*}{$ \frac{\tau_k}{f_{k,l}^{T_{\z}}-f(\x_l, \y_{k,l}^{T_{\y}} )}{\left(\frac{\partial f(\x_l,\y_{k,l}^{T_{\y}})}{\partial \x} - \frac{\partial f(\x_l,\z_{k,l}^{T_{\z}})}{\partial \x} \right)}$}}&\multirow{3}{*}{$O(c(T_{\z}+T_{\y}))$}&\multirow{3}{*}{$O(m+n)$}\\
				&&&\\     
				 &&&\\                                                               
				\bottomrule
			\end{tabular}
		\end{small}
	\end{center}

\end{table*}

\section{ Complexity Analysis}\label{complexity}
	In this part, we compare the time and space complexity of Algorithms~\ref{alg:ours} in computing the direction for updating variable $\x$ with EGBMs (i.e.,FHG and RHG) and IGBMs (i.e.,CG and Neumann). Our complexity analysis is based on the basic step. 
	We suppose the gradients of $F$, $f$, Hessian-vector product $\frac{\partial^2  f}{\partial \y^2} \q$ and Jacobian-vector product $\frac{\partial^2 f}{\partial \x \partial \y} \q$ can be evaluated in time $c= c(m,n)$ for any vector $\q \in \mathbb{R}^n$. 
	For all existing methods, we assume the optimal solution of the LL problem is calculated by a $T$-step gradient descent, and the transition function in FHG and RHG, also the $\Phi$ in Eq.~\eqref{dynamic}, is gradient descent. 
	In practice, it has been shown in~\cite{rajeswaran2019meta} that the time and space cost for computing $ \frac{\partial^2 f}{\partial \y^2} \q$ and $\frac{\partial^2 f}{\partial \x \partial \y}\q$ via AD is no more than a (universal) constant (e.g., usually 2-5 times) order over the cost for computing gradient of $F$ and $f$.

	\textbf{EGBMs.} As discussed in~\cite{franceschi2017forward,shaban2019truncated}, FHG for forward calculation of Hessian-matrix product can be evaluated in time $O(cmT)$ and  space $O(mn)$, and RHG for reverse calculation of Hessian- and Jacobian-vector products can be evaluated in time $O (cT)$  and  space $O (m+nT)$. The time complexity for EGBMs to calculate the UL gradient is proportional to the number of iterations $T$ of the LL problem, and thus EGBMs take a large amount of time to ensure convergence.

	\textbf{IGBMs.} After implementing a $T$-step gradient descent for the LL problem, IGBMs approximate matrix inverse by solving $J$ steps of a linear system. As each step contains a Hessian-vector product, this part required $O(cJ) $ time and $O(m+n)$ space. So IGBMs run in time $O (c(T+J)) $ and space $O (m+n)$. 
	The iteration number $J$ always relies on the properties of the Hessian-matrix, and thus in some cases, $J$ can be much larger than $T$.

	\textbf{BVFIM.} In our algorithm, $T_{\z}$ steps of gradient descent in Eq.~\eqref{update_z} for the LL problem solution $\z^{T_{\z}}$ requires time  $O(cT_{\z})$ and space $O(n)$. Then $T_{\y}$ gradient descent steps in Eq.~\eqref{update_y} are used calculate $\y^{T_{\y}}$, which requires time $O(cT_{\y})$ and space $O(n)$.
	Lastly, the direction can be obtained according to the formula given in Eq.~\eqref{a_g_vphi} by several computations of the gradient $\frac{\partial f}{\partial \x}$ and $\frac{\partial F}{\partial \x}$ without any intermediate update, which requires time $O(c)$ and space $O(m+n)$, so BVFIM runs in time $O(c(T_{\z}+T_{\y}))$ and space $O(m+n)$, respectively. Note that although BVFIM has similar complexity to existing methods, 
	it does not need any computation of Hessian- or Jacobian-vector products and all the complexity comes from calculating the gradients of $F$ and $f$. 
	And because calculating the gradient is much faster than calculating Hessian- and Jacobian-vector products (even by AD) when the dimension of LL variable is large, the proposed BVFIM is more efficient than existing methods when the LL problem is high-dimensional, for example, in applications with a large-scale network. We will illustrate this advantage in the next section.
	\begin{figure*}  
		\centering  
		\includegraphics[height=3cm,width=4cm]{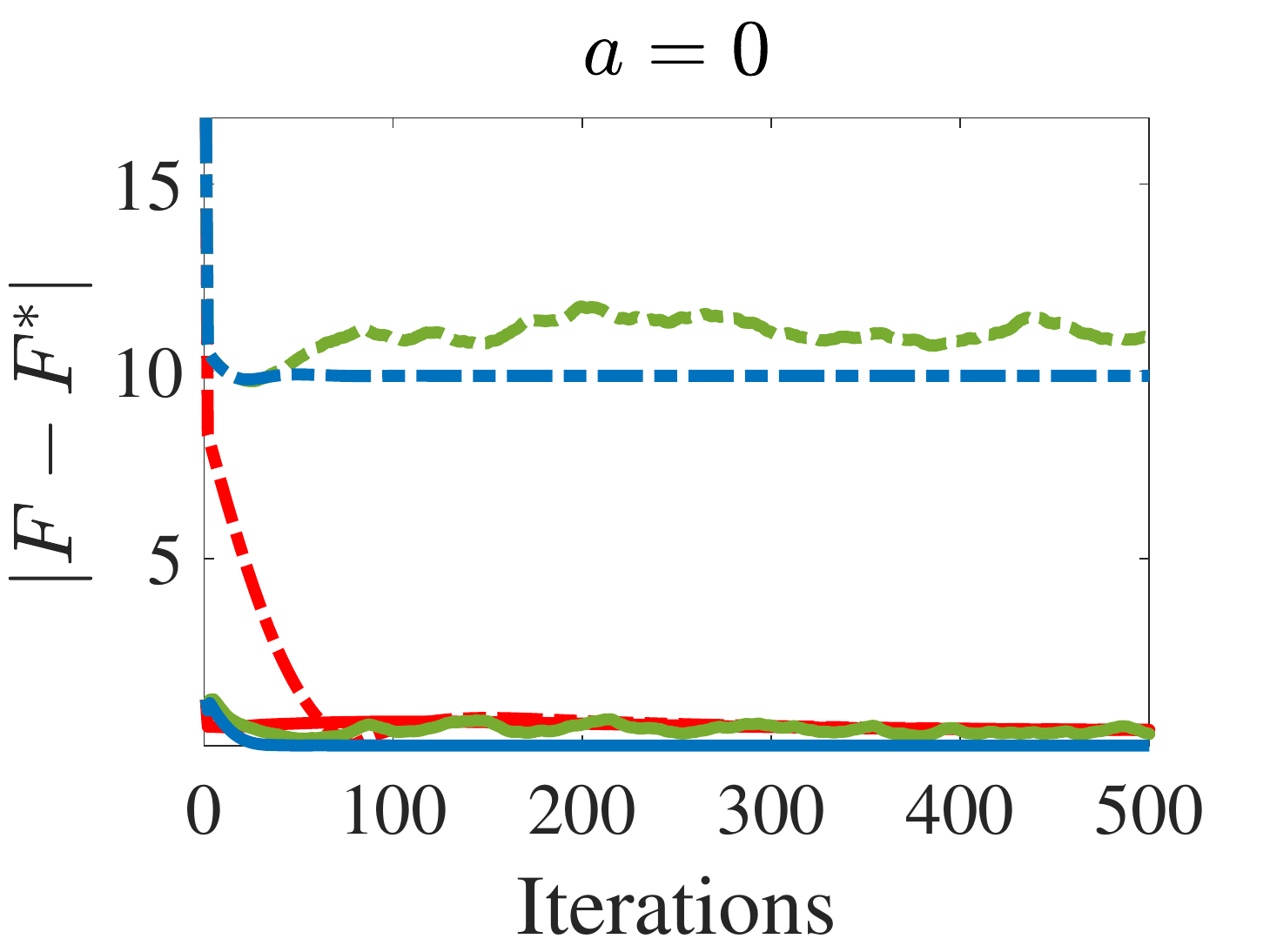}  
		\includegraphics[height=3cm,width=4cm]{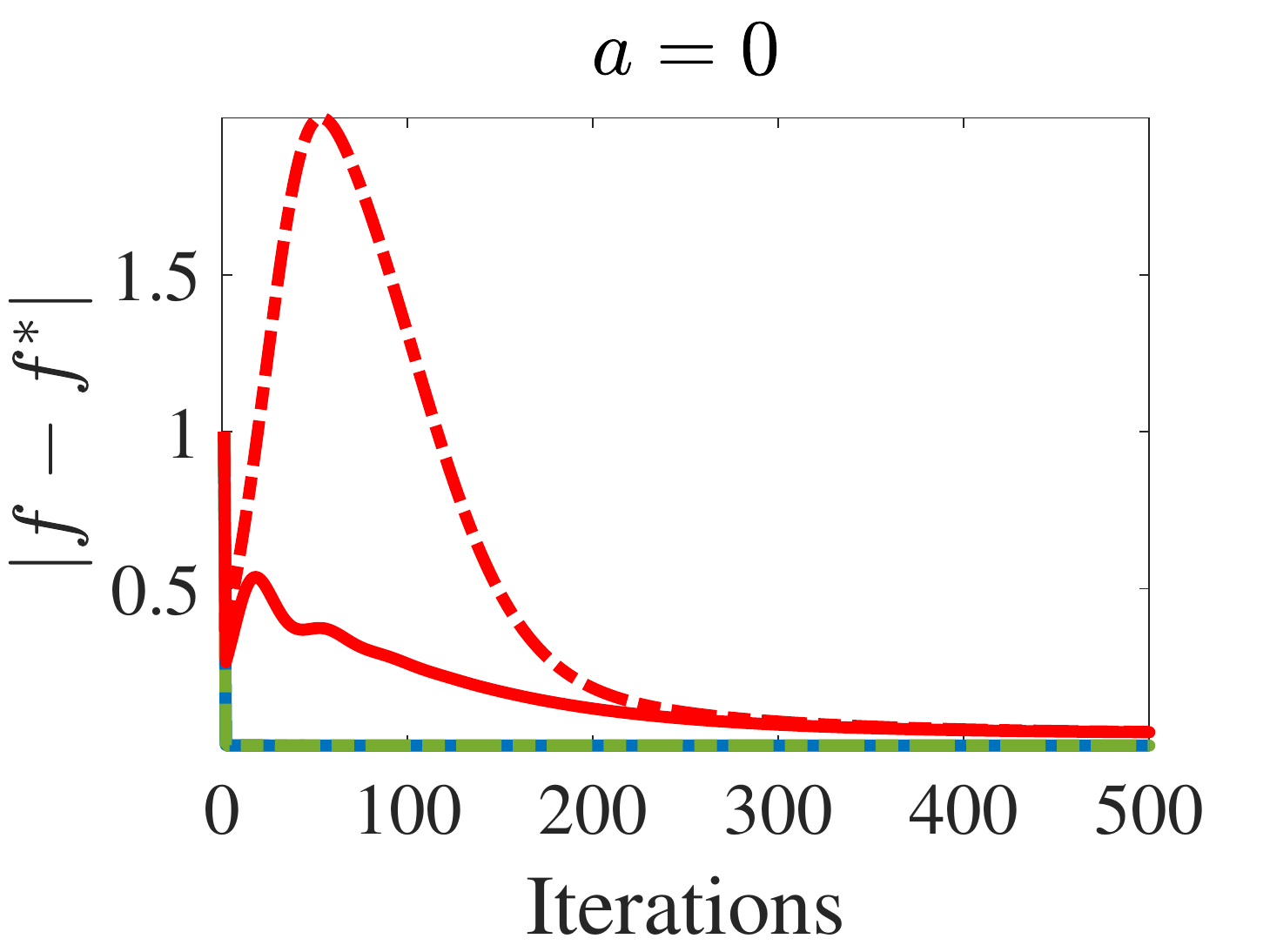}  
		\includegraphics[height=3cm,width=4cm]{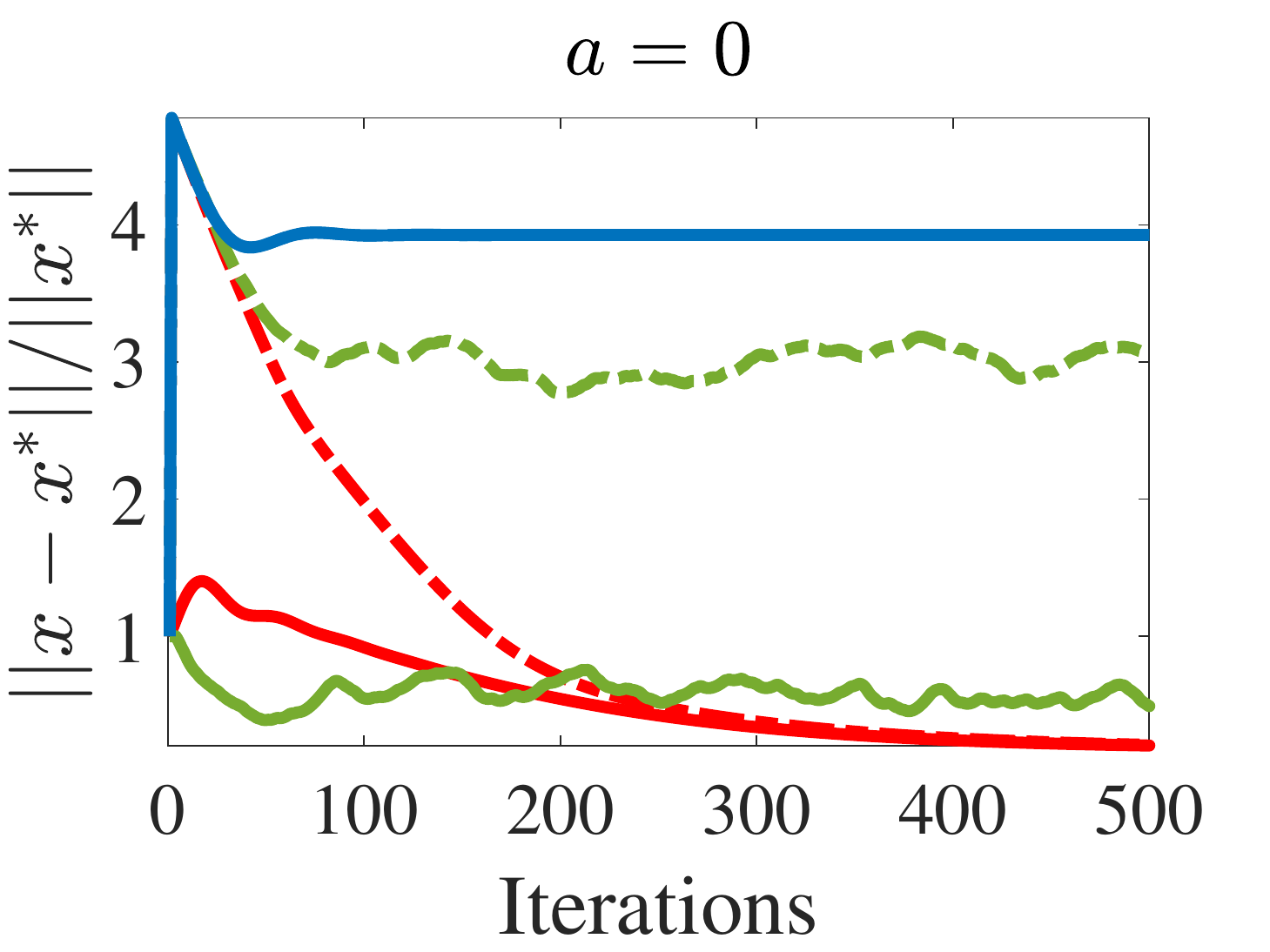}  
		\includegraphics[height=3cm,width=4cm]{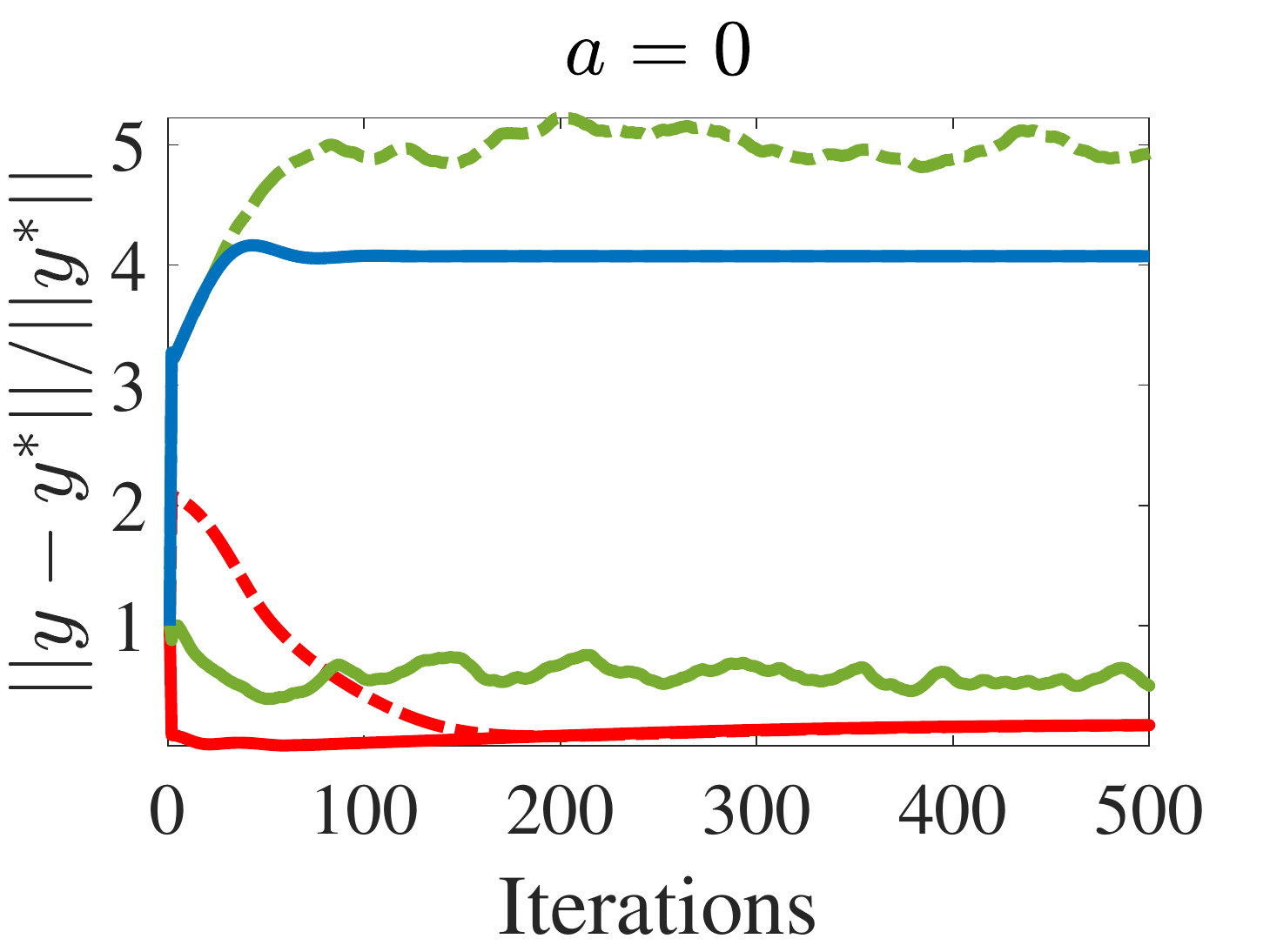}  
		\includegraphics[height=3cm,width=4cm]{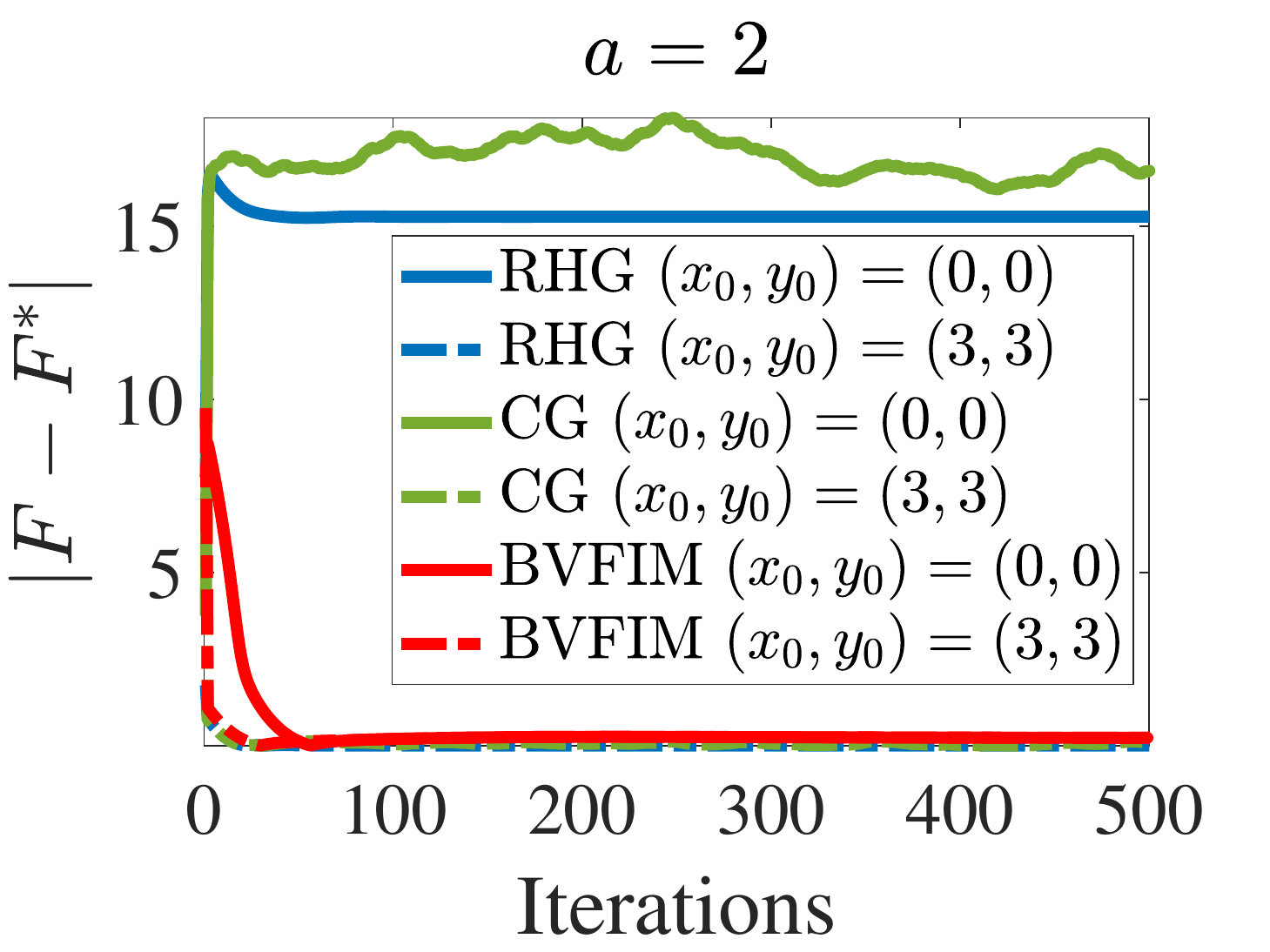}  
		\includegraphics[height=3cm,width=4cm]{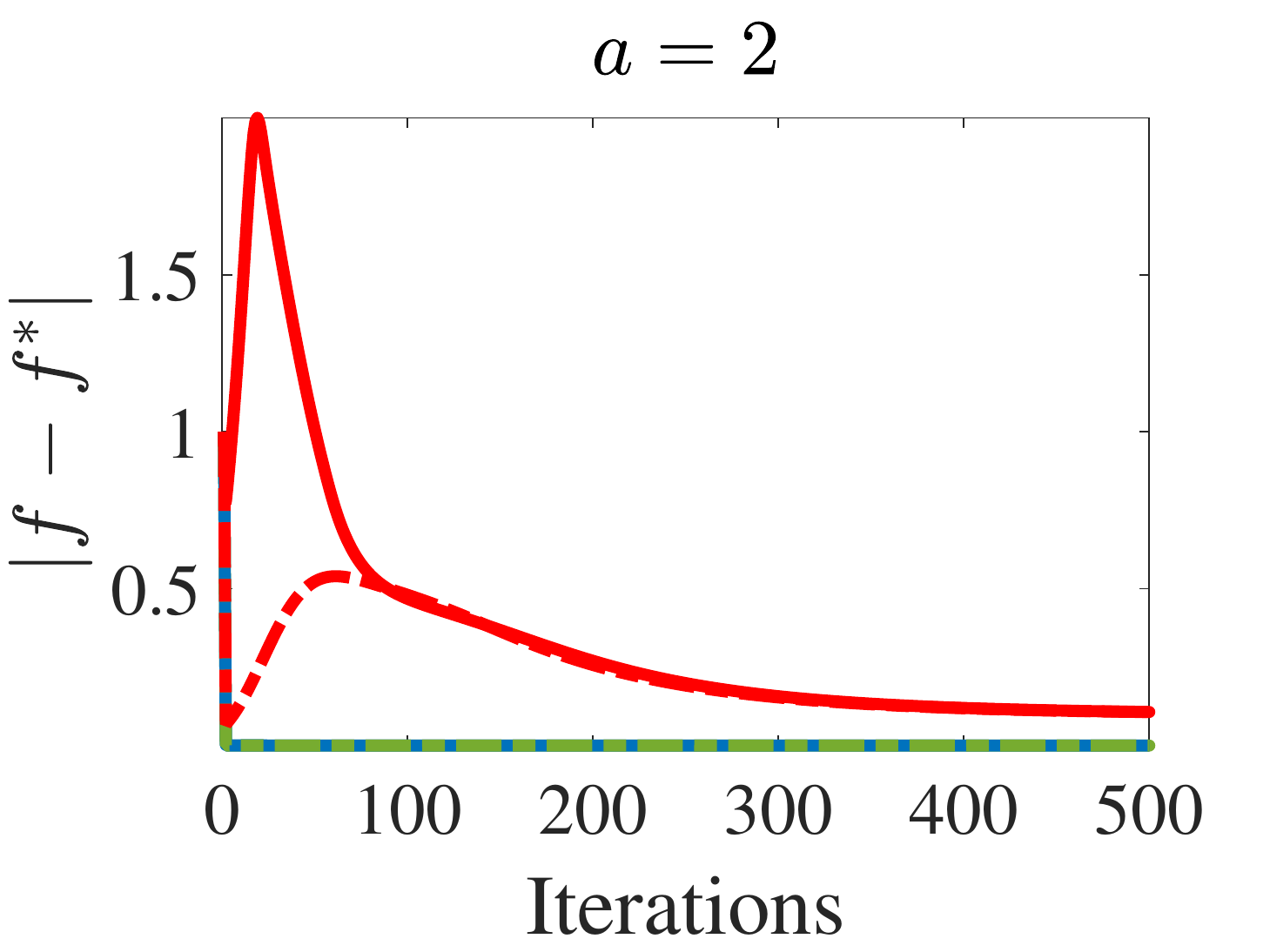}  
		\includegraphics[height=3cm,width=4cm]{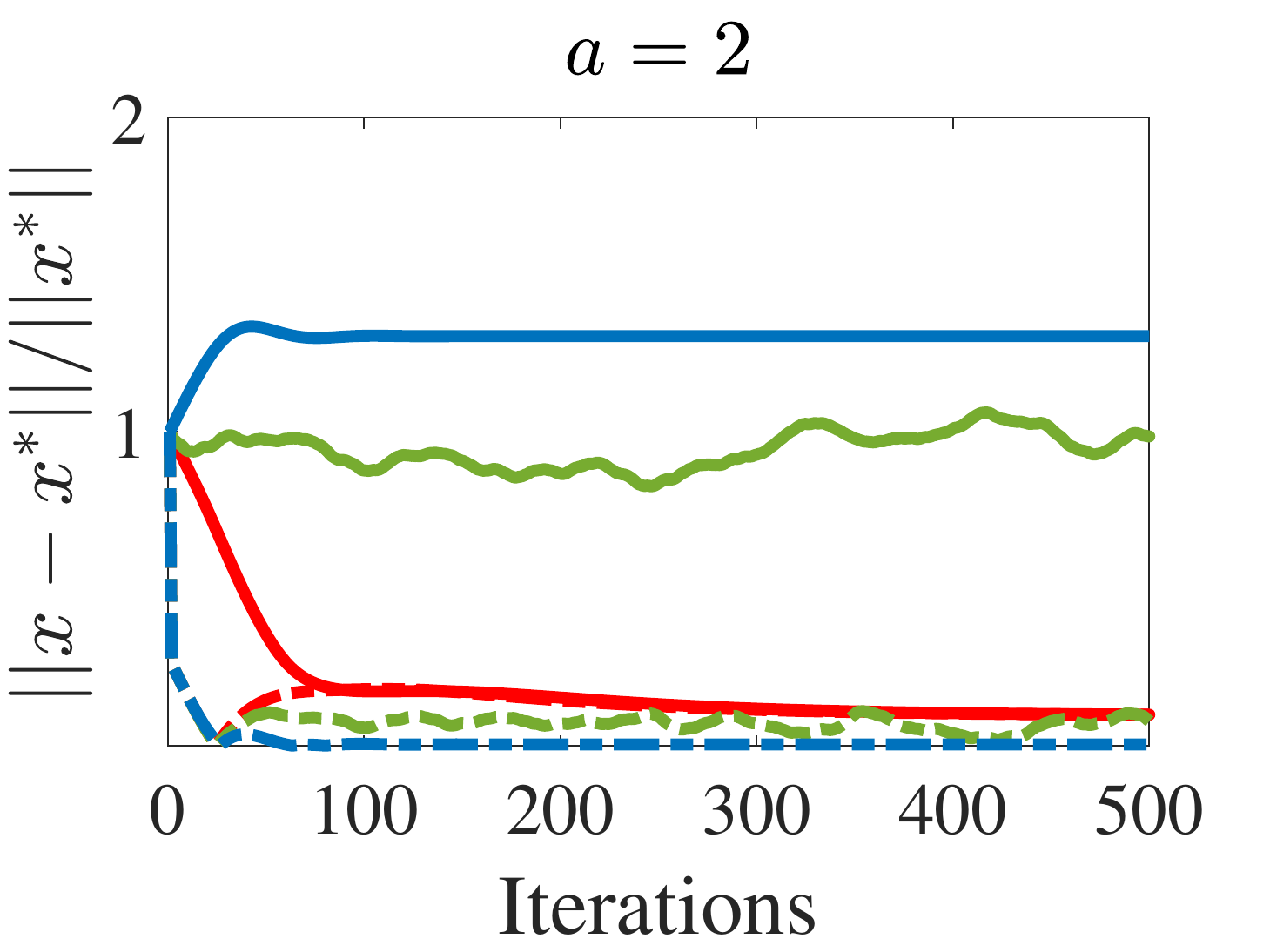}  
		\includegraphics[height=3cm,width=4cm]{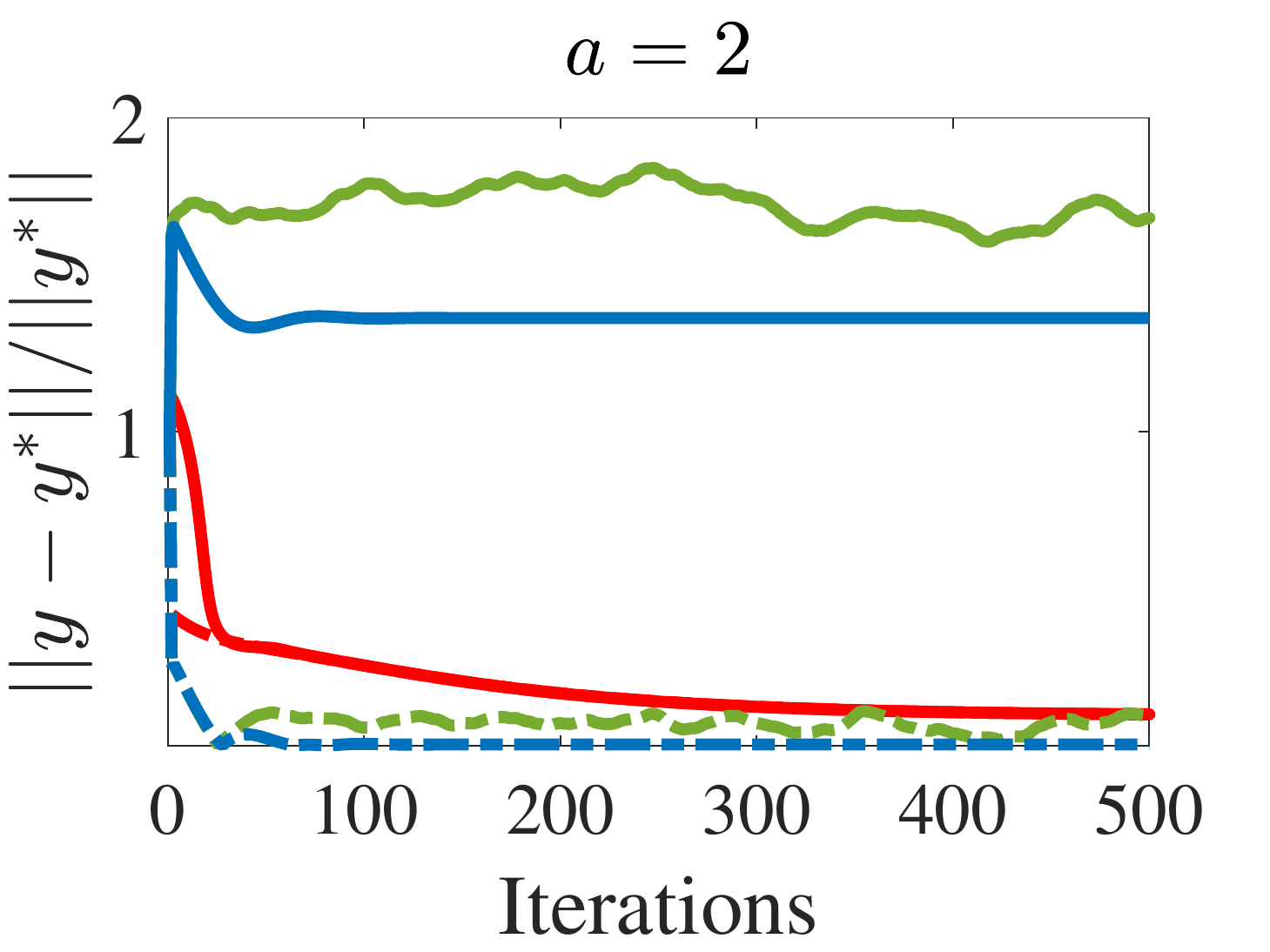}  
		\caption{Comparison of the existing explicit method (RHG) and implicit method (CG) with BVFIM on different metrics and different initial points. We set $a = 0$ in the first row and $a = 2$ in the second row. The solid line represents $(x_0,y_0)=(0,0)$, and the dotted line represents $(x_0,y_0)=(3,3)$.  The legend is only plotted in the first subfigure of the second line.}  
		\label{1}  

	\end{figure*} 
\section{Experimental Results}\label{experiment}
	 In this section, we present numerical simulations with the
	proposed methods and demonstrate application in hyper-parameter optimization for non-convex LL problems.

	\subsection{Toy Examples}
	In this section, our goal is to use a simple numerical sample to verify the superiority of BVFIM over existing methods in the case of non-convex LL problems. We use the following numerical examples \footnote{Since variables are one-dimensional real numbers in numerical experiments, we do not use bold letters which represent vectors.	}
	\begin{equation}
	\min_{ x \in \mathbb{R}, y \in \mathbb{R}} \Vert x-a\Vert^2+\Vert y-a\Vert^2, \text{\ s.t.\ }\; y \in \underset{ y \in \mathbb{R}}{\mathrm{argmin}}\; \sin( x+ y),\label{non-convexExperiment}
	\end{equation}
	where $a$ is an adjustable parameter which will be set as $a = 0$ and $2$ in following experiments. 
	By simple calculation, we know that the optimal solution of Eq.~\eqref{non-convexExperiment} is $ x^*_{a=0}= y^*_{a=0}=-\pi/4$, $F^*_{a=0}=\pi^2/8$ and $ x^*_{a=2}= y^*_{a=2}=3\pi/4$, $F^*_{a=3}=3\pi^2/16-9\pi/2+9$. 
	This example satisfies all the assumptions of BVFIM, but not LLS and LLC assumptions in ~\cite{pedregosa2016hyperparameter,rajeswaran2019meta,lorraine2020optimizing}, which makes it be a good example to validate the advantages of BVFIM.

	In Figure~\ref{1}, we show the comparison of BVFIM with the explicit method RHG and the implicit method CG for solving different bi-level optimization problems at different initial points. It can be observed that RHG and CG converge to the optimal point only when the initial point is near the optimal point, but different initial points only slightly affect the convergence of BVFIM. From the second column in Figure~\ref{1}, we can see that the effectiveness of the BVFIM method in solving non-convex problems comes at the expense of slowing down the convergence rate of the LL problem $f$.

	\begin{figure} 
		\centering  
		\includegraphics[height=3cm,width=3.9cm]{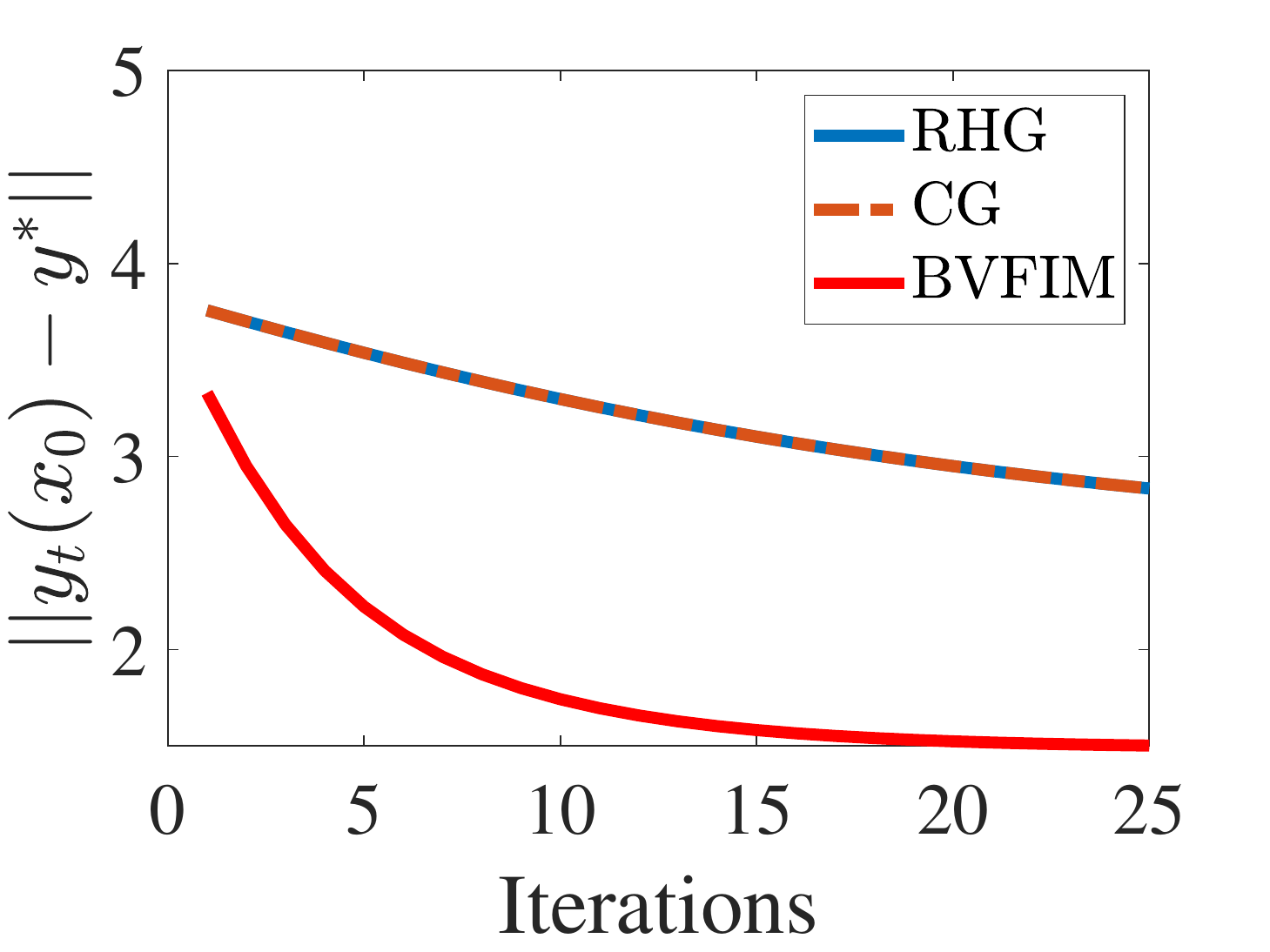}  
		\includegraphics[height=3cm,width=4cm]{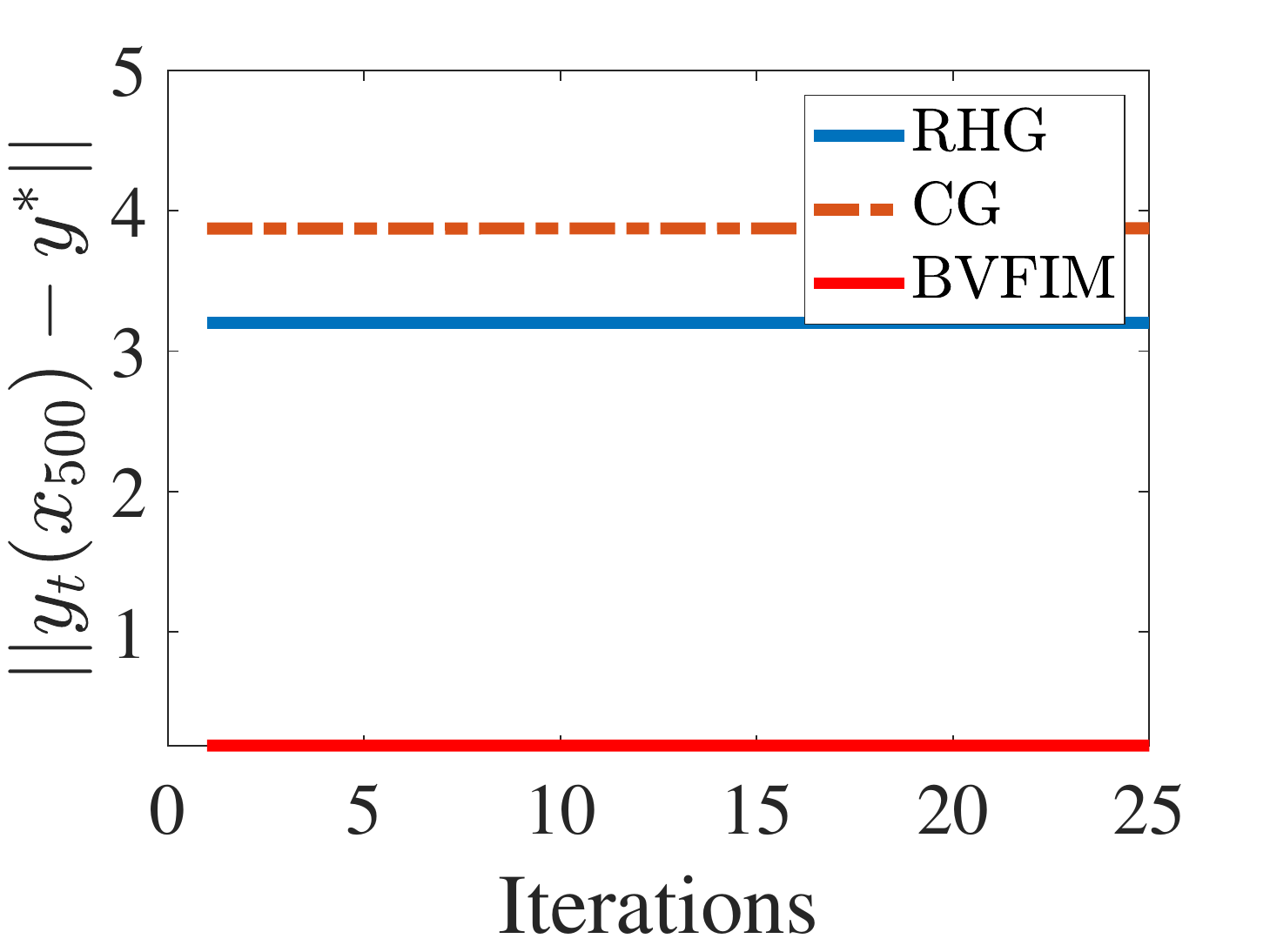}  
		\includegraphics[height=3cm,width=4cm]{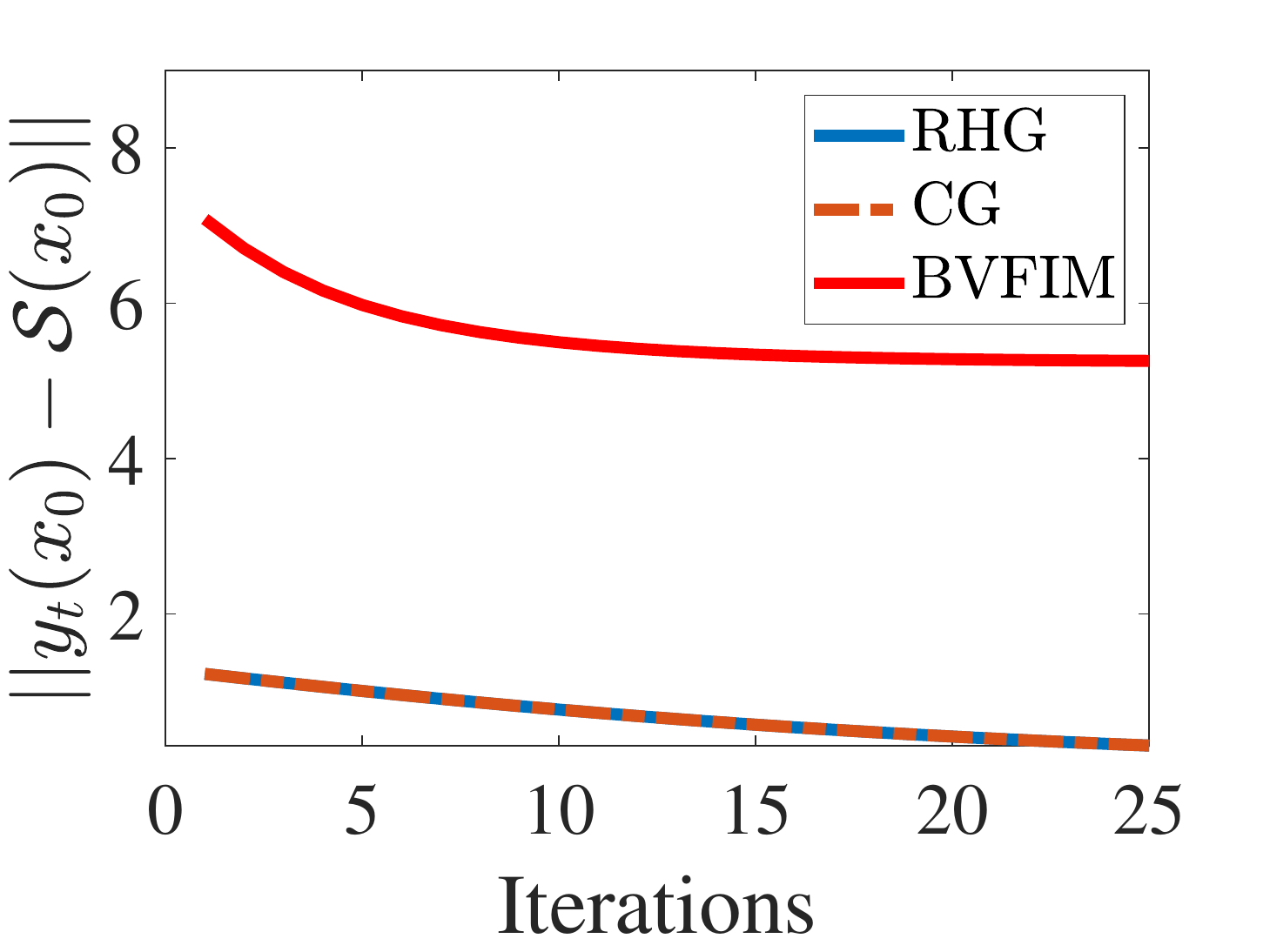}  
		\includegraphics[height=3cm,width=4cm]{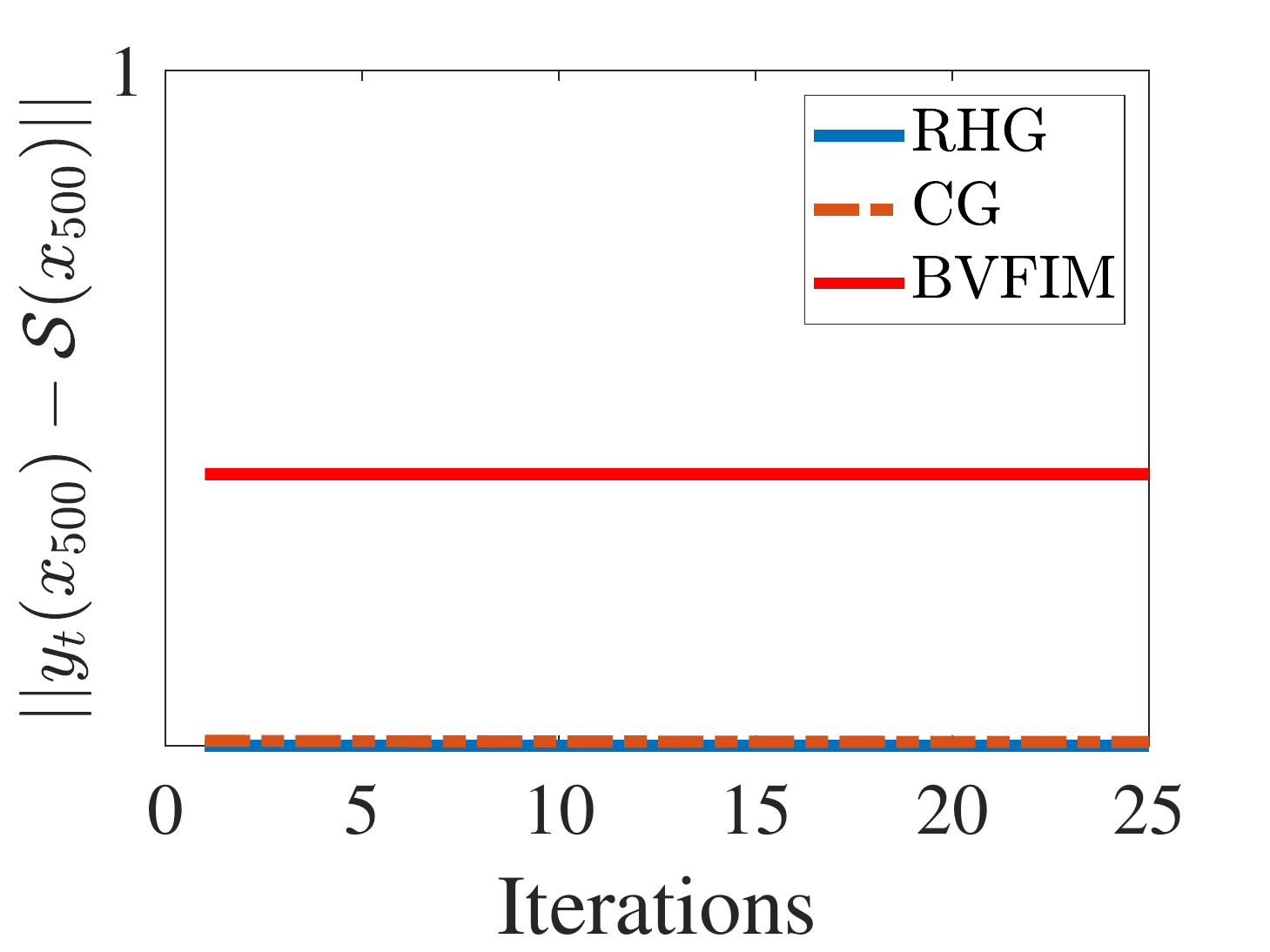}  
		\caption{LL iterative curves with explicit and implicit methods (RHG and CG) in $ x_0$ and $ x_{500}$. $\S ( x) $ is the solution set of the LL problem with given $ x$, and $ y^*$ is the global optimal solution. We set $a=0$ and $( x_0, y_0)=(3,3)$.}  
		\label{2}  

	\end{figure} 
	
	\begin{figure} 
		\centering  
		\includegraphics[height=3cm,width=4cm]{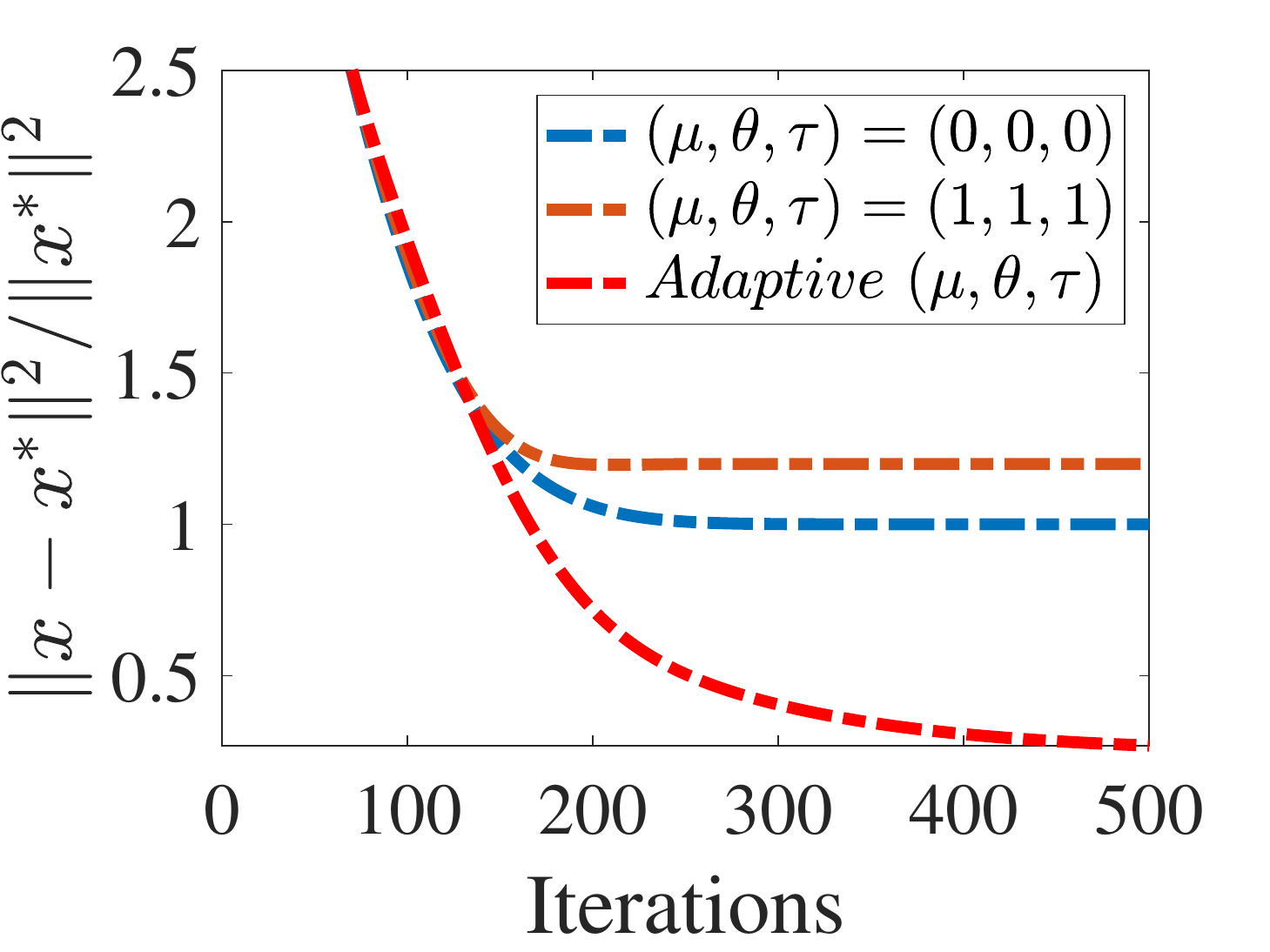}  
		\includegraphics[height=3cm,width=4cm]{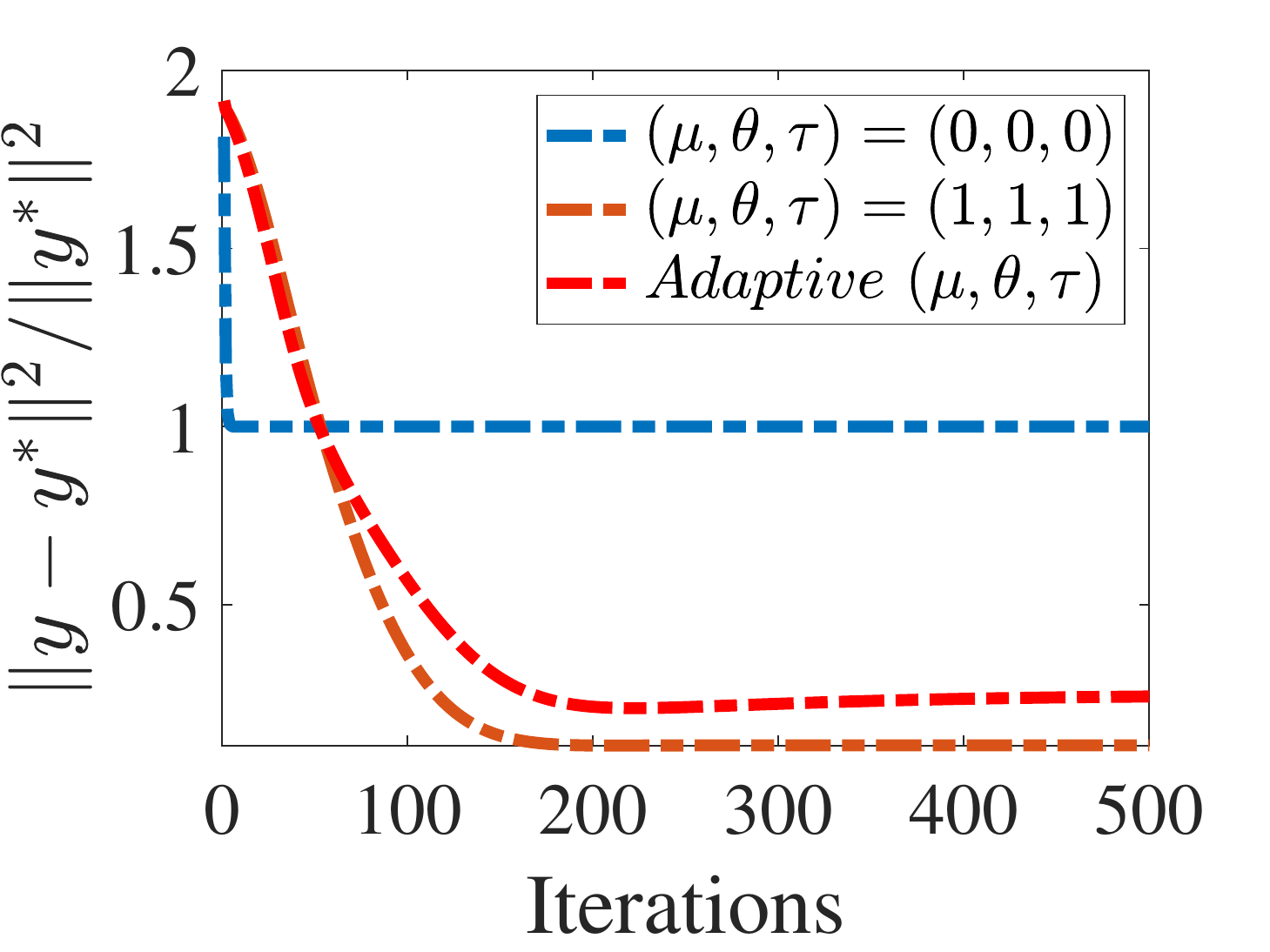}  
		\caption{Convergence curves under different $(\mu_k,\theta_k,\tau_k)$ settings. The legend is only plotted in the second subfigure. We set $a=0$ and $( x_0, y_0)=(3,3)$. The adaptive $(\mu,\theta,\tau)$ means to follow the default setting which can be found in Supplementary Materials. } %
		\label{3}  
	\end{figure} 
	To further illustrate the differences between the BVFIM methods and existing methods on the LL problem, we plot the convergence behavior of the LL iteration (i.e., $||y_{t}(x_0)-y^*||$ and $||y_{t}(x_{500})-\mathcal{S}(x_{500})||$) with given UL variable $ x$ in Figure~\ref{2}. Note that $ x_0$ and $ x_{500}$ refer to the $ x$ at the 0th and 500th iteration of the UL problem, respectively. From the first row in Figure~\ref{2}, it can be observed that RHG and CG converge quickly to the solution set of the LL problem, while BVFIM does not because $ y_t( x)$ in BVFIM is solving problem Eq.~\eqref{BVFIM}. However, the second row of Figure~\ref{2} shows that BVFIM can converge to the global optimal solution while RHG and CG cannot.
	
			\begin{table*}[t]
		\caption{Comparison of the results of existing methods for solving data hyper-cleaning tasks on MNIST, FashionMNIST and CIFAR10. The F1 score denotes the harmonic mean of the precision and recall.}
		\begin{center}
			\begin{small}
				\begin{tabular}{lccccccccc}
					\toprule		
					\multirow{2}{*}{Method}&\multicolumn{3}{c}{MNIST}&\multicolumn{3}{c}{FashionMNIST}&\multicolumn{3}{c}{CIFAR10}\\
					&Acc.&F1 score&Time(s)&Acc.&F1 score&Time(s)&Acc.&F1 score&Time/s\\
					\midrule
					RHG&87.90&89.36&0.4131&81.91&87.12&0.4589&34.95&68.27&1.3374\\
					TRHG&88.57&89.77&0.2623&81.85&86.76&0.2840&35.42&68.06&0.8409\\
					CG&89.19&85.96&0.1799&83.15&85.13&0.2041&34.16&69.10&0.4796\\
					Neumann&87.54&89.58&0.1723&81.37&87.28&0.1958&33.45&68.87&0.4694\\
					BDA&87.15&90.38&0.6694&79.97&88.24&0.8571&36.41&67.33&1.4869\\
					BVFIM&\textbf{90.41}&\textbf{91.19}&\textbf{0.1480}&\textbf{84.31}&\textbf{88.35}&\textbf{0.1612}&\textbf{38.19}&\textbf{69.55}&\textbf{0.4092}\\
					\bottomrule	
				\end{tabular}
			\end{small}
		\end{center}
	\end{table*}	
	Figure~\ref{3} shows the convergence behavior of BVFIM with different regularization parameter settings. It can be seen that by setting the regularization parameter to 0, our method degenerates to optimize UL and LL problems separately without regard to their relevance, thus it is difficult to obtain an appropriate feasible solution. Choosing a fixed nonzero regularization parameter can improve convergence on $ y$ and achieve the same convergence behavior on $ x$ as adaptive $\left( \mu,\theta,\tau\right) $ at the early stage. However, due to the influence of the regular term, it does not converge to the optimal solution on $ x$. When the regularization term takes a decreasing value, the best results are obtained.
	\begin{figure}  
		\centering  
		\includegraphics[height=3cm,width=4cm]{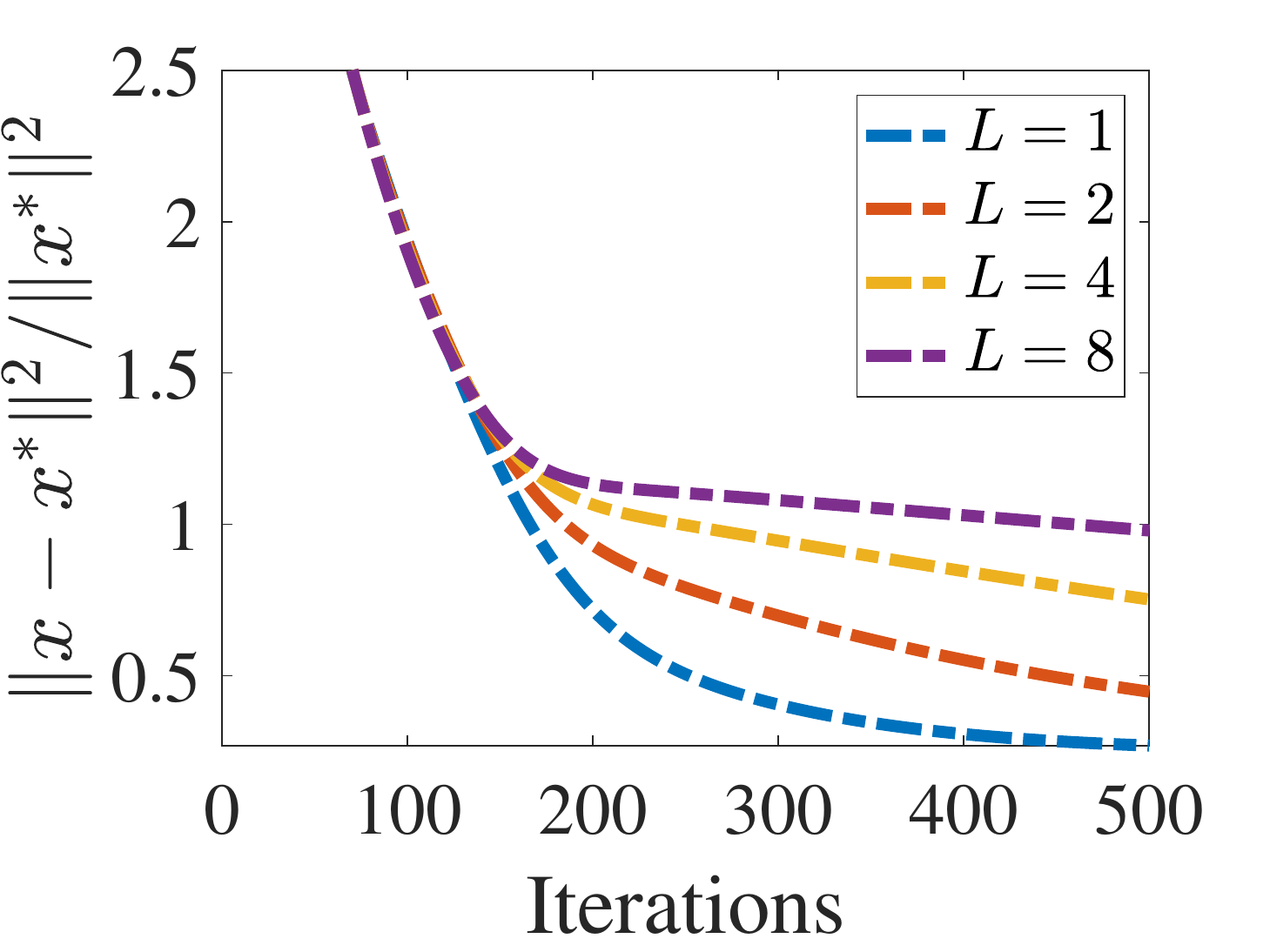}  
		\includegraphics[height=3cm,width=4cm]{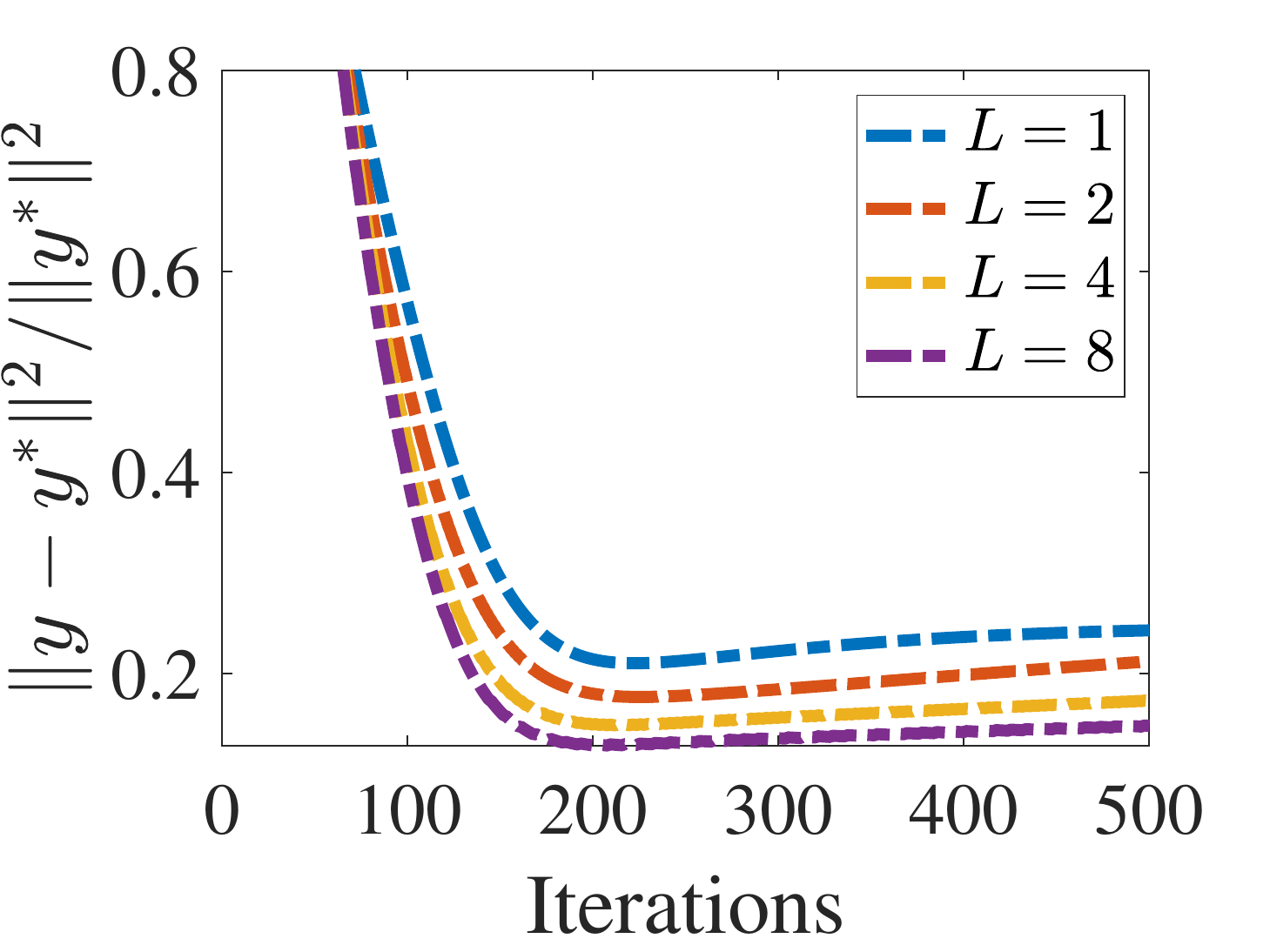}  
		\caption{Comparison of convergence rates under different $L$ settings. We set $a=0$ and $( x_0, y_0)=(3,3)$. The iterations are the number of gradient descents of $ x$, which is $K\times L$.}  
		\label{4} 
	\end{figure}

	Figure~\ref{4} shows the effect of different iteration number $L$ on convergence behavior.
	We can find that smaller $L$ can accelerate the convergence of variable $ x$, while larger $L$ can improve the convergence of variable $ y$.
	These preliminary results indicate that in situations where optimal solution $ x^*$ is needed, BVFIM with small $L$ (e.g. $L=1$) is a better choice.
	

	\subsection{Hyper-parameter Optimization}

	In this section, we use a special case of hyper-parameter optimization, called data hyper-cleaning, to evaluate the performance of our method when the LL problem is non-convex.
	Assuming that some of the labels in our datasets are contaminated, the goal of data hyper-cleaning is to reduce the impact of incorrect samples by adding hyper-parameters to the contaminated samples.
	
	In the experiment, we set $\y \in \mathbb{R}^{10\times301}\times\mathbb{R}^{300\times d}$ as the parameters of 2-layer linear network classifier where $d$ is the data dimension and  $\x\in \mathbb{R}^{|\D_{\mathtt{tr}}|}$ as the weight of each sample in the training set. Therefore, the LL problem is to learn a classifier $\y$ on cross-entropy
	training loss weighted with given $\x$
	\begin{equation}
	f(\x,\y)=\sum_{(u_{i},v_{i})\in \D_{\mathtt{tr}}}\sigma(\x_i)\mathtt{CE}(\y,u_{i},v_{i}),
	\end{equation}
	where  $(u_{i},v_{i})$ are the training examples, $\sigma(\x)$ denotes the sigmoid function on $\x$ to constrain the weights in the range $[0, 1]$ and $\mathtt{CE}$ is the cross entropy loss. The UL problem is to find a weight $\x$ which can reduce the cross-entropy loss of $\y$ on a cleanly labeled validation set 
	\begin{equation}
	F(\x,\y)=\sum_{(u_{i},v_{i})\in \D_{\mathtt{val}}}\mathtt{CE}(\y,u_{i},v_{i}).
	\end{equation}
	Our experiment is based on MNIST, FashionMNIST ~~\cite{xiao2017fashion} and CIFAR10 datasets.
	For each dataset, we randomly select 5000 samples as the training set $\D_{\mathtt{tr}}$, 5000 samples as the validation set  $\D_{\mathtt{val}}$, and 10000 samples as the test set $\D_{\mathtt{test}}$. Then, we contaminated half of the labels in $\D_{\mathtt{tr}}$. 
	\begin{figure}  
		\centering  
		\includegraphics[height=3cm,width=4cm]{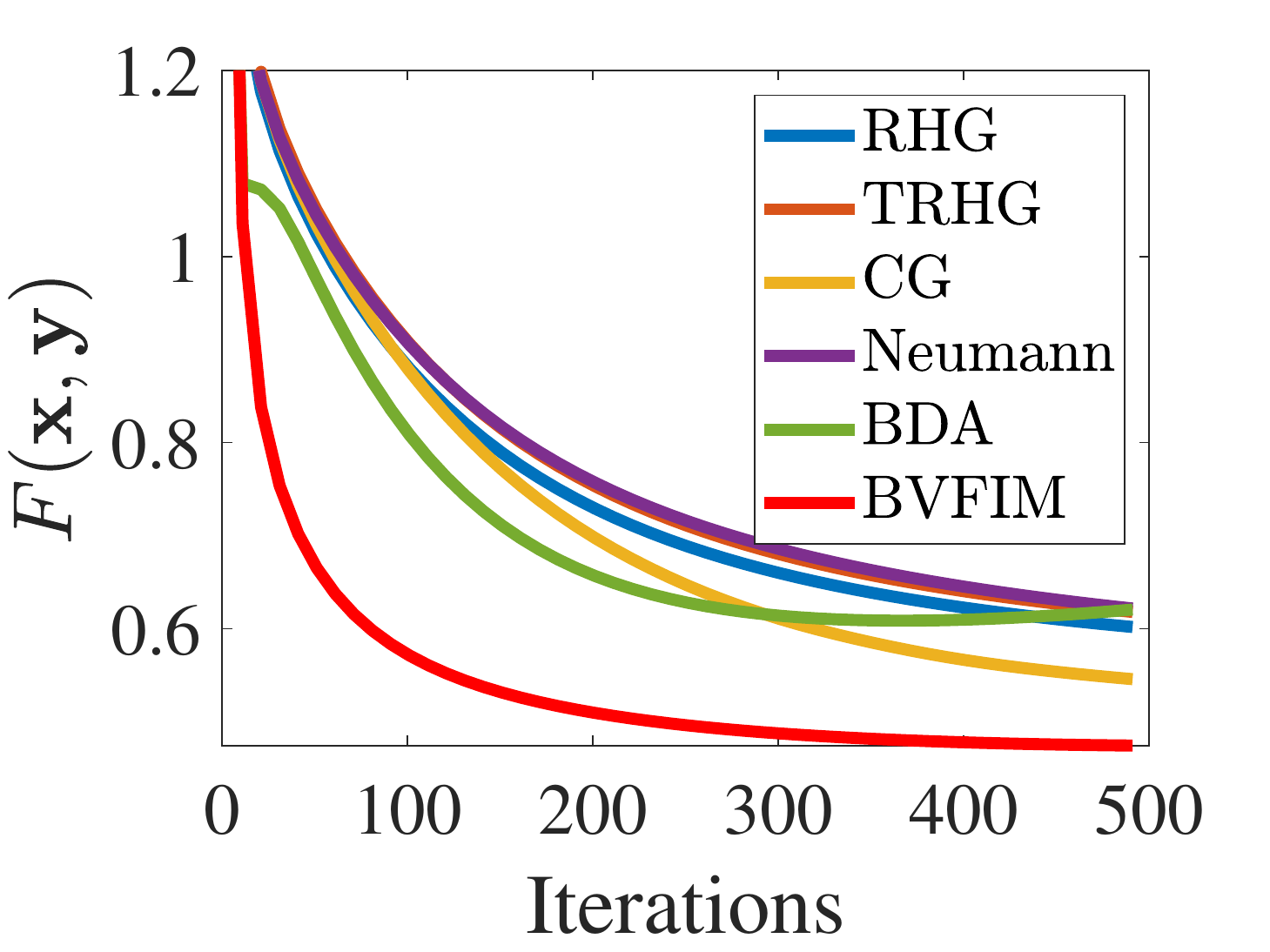}  
		\includegraphics[height=3cm,width=4cm]{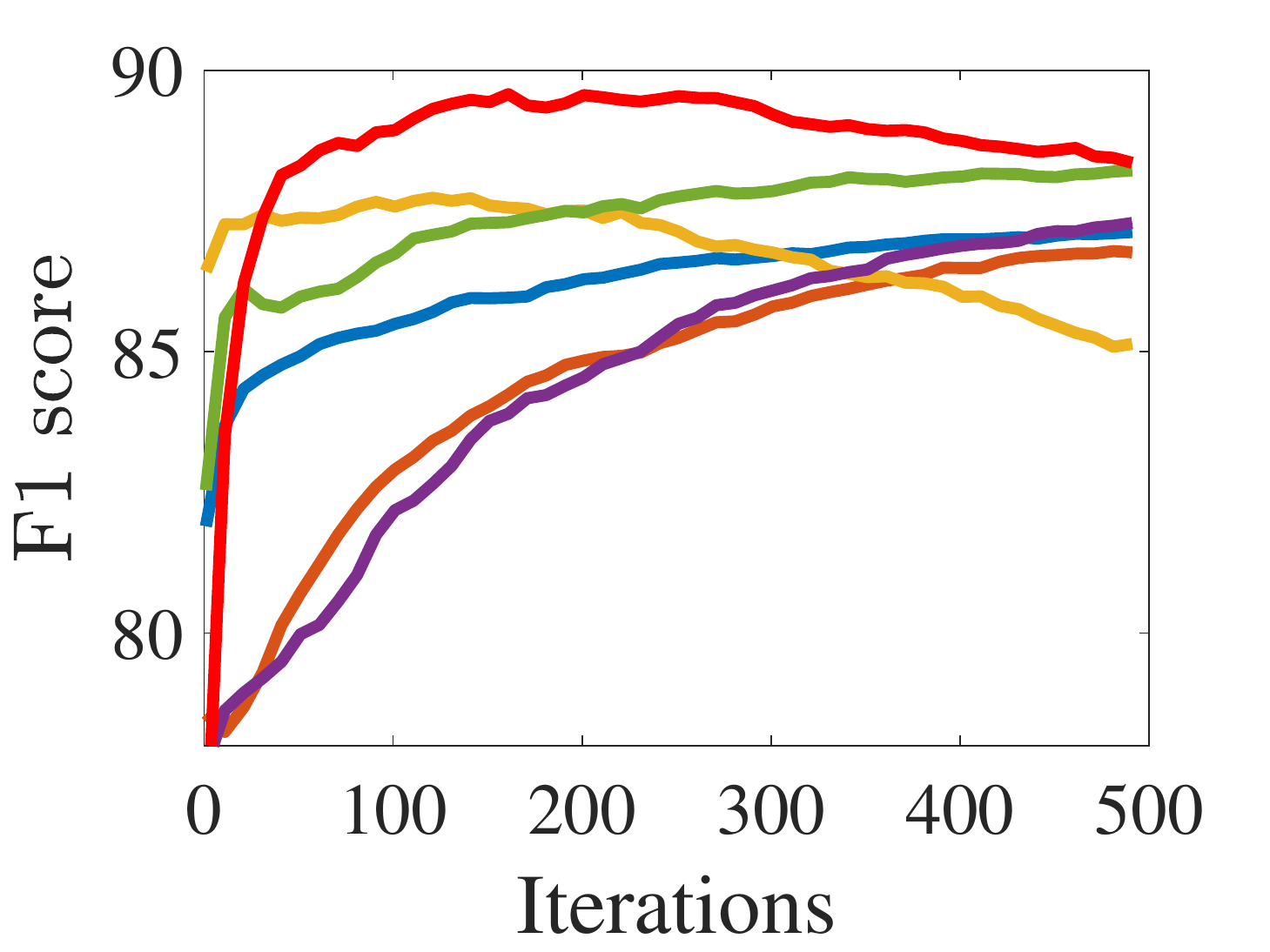}  
		\caption{$F(\x,\y)$ and F1 score between existing methods and BVFIM. The curves are based on the FashionMNIST experiment in Table 3. The legend is only plotted in the first subfigure.}  
		\label{5}  
	\end{figure}

	Table 3 shows the results and calculation time of the compared methods and BVFIM on three different datasets. It can be seen that BVFIM achieves the most competitive results on all three datasets. 
	Furthermore, BVFIM is faster than EGBMs and IGBMs, and this advantage is more evident in CIFAR10 datasets with larger LL parameters, consistent with the complexity analysis given above. The UL objective and F1 scores of BVFIM and compared methods on the FashionMNIST dataset are plotted in Figure~\ref{5}.

	\section{Conclusions}
	In this paper, we present a new bi-level optimization algorithm BVFIM. By penalizing the regularized LL problem into the UL problem, we approximate the BLO problem by a series of unconstrained smooth SLO sub-problems. We prove the convergence without requiring any convexity assumption on either the UL or the LL objectives. Extensive evaluations show the superiority of BVFIM for different applications.

	\section*{Acknowledgements}
	
This work is partially supported by the National Natural Science Foundation of China (Nos. 61922019, 61733002 and 11971220), LiaoNing Revitalization Talents Program (No. XLYC1807088),  the Fundamental Research Funds for the Central Universities, Shenzhen Science and Technology Program (No. RCYX20200714114700072), and Guangdong Basic and Applied Basic Research Foundation (No. 2019A1515011152).

	\bibliography{output}
	\bibliographystyle{icml2021}

\appendix
	
\section*{Appendix}
	This supplementary material is organized as follows. In Section~\ref{sectionA} and Section~\ref{sectionB}, we present the detailed proof of Section 3  and subsection 4, respectively.
    Section \ref{sectionC} presents the further details of the experiments in
	Section \ref{experiment}.  Section \ref{sectionD} presents some additional results.
	\section{Proof of Section \ref{secBVFIM}}\label{sectionA}
	\subsection{Proof of Proposition \ref{gradient}}
	\begin{proof}
		Since ${\mathrm{argmin}}_{\y \in \mathbb{R}^n}~f(\x,\y)  + \frac{\mu_1}{2}\|\y\|^2 + \mu_2$ is a singleton, it follows from \cite{AlexanderShapiro2011Perturbation}[Theorem 4.13, Remark 4.14] that
		\begin{equation}
		\begin{aligned}
		\frac{\partial  f_{\mu}^*(\x)}{\partial\x} &= \left. \frac{\partial \left( f(\x,\y) + \frac{\mu_{1}}{2}\|\y\|^2 + \mu_{2}   \right)}{\partial \x}  \right|_{\y = \z_{\mu}^*(\x)} \\
		&= \frac{\partial f(\x, \z_{\mu}^*(\x))}{\partial\x}.
		\end{aligned}
		\end{equation}
		As ${\mathrm{argmin}}_{\y\in\mathbb{R}^n}~F(\x,\y) + \frac{\theta}{2}\|\y\|^2 - \tau\ln\left(f_{\mu}^*(\x) - f(\x,\y)\right)$ is a singleton, ~\cite{AlexanderShapiro2011Perturbation}[Theorem 4.13, Remark 4.14] shows that
		\begin{equation}
			\begin{aligned}\small
			&\frac{\partial \varphi_{\mu,\theta,\tau}(\x)}{\partial\x} =  \\
			& \frac{\partial\left( F(\x,\y) + \frac{\theta}{2}\|\y\|^2 - \tau\ln\left(f_{\mu}^*(\x) - f(\x,\y)\right) \right)}{\partial \x}\left. \right|_{\y = \y_{\mu,\theta,\tau}^*(\x)}\\
			&= \frac{\partial F(\x,\y_{\mu,\theta,\tau}^*(\x))}{\partial\x}	 + \frac{\tau \left( \frac{\partial f(\x,\y_{\mu,\theta,\tau}^*(\x))}{\partial\x} - \frac{\partial f_{\mu}^*(\x) }{\partial\x}\right)}{f_{\mu}^*(\x) - f(\x,\y_{\mu,\theta,\tau}^*(\x))} \\
			& =  \frac{\partial F(\x,\y_{\mu,\theta,\tau}^*(\x))}{\partial\x}+ \frac{\tau\left( \frac{\partial f(\x,\y_{\mu,\theta,\tau}^*(\x))}{\partial\x} - \frac{\partial f(\x,\z_{\mu}^*(\x))}{\partial\x} \right)}{f_{\mu}^*(\x) - f(\x,\y_{\mu,\theta,\tau}^*(\x))}.
			\end{aligned}
		\end{equation}
	\end{proof}
	\section{Proof of Section \ref{secTheo}}\label{sectionB}
	\subsection{Proof of Lemma \ref{lem2}}
	\begin{proof}
		For any $\epsilon > 0$, there exists $\bar{\y} \in \mathbb{R}^n$ such that $f(\bar{\x},\bar{\y}) \le f^*(\bar{\x}) + \epsilon$. And as $\mu_k \rightarrow 0$, we can find $k_1 > 0$ such that
	\begin{equation}	
		f(\bar{\x},\bar{\y}) + \frac{\mu_{k,1}}{2}\|\bar{\y}\|^2 + \mu_{k,2} \le f^*(\x) + 2\epsilon
		\end{equation}
		for all $k \ge k_1$. Next, as $\{\x_k\}$ converging to $\bar{\x}$, it follows from the continuity of $f(\x,\y)$ that there exists $k_2 \ge k_1$ such that $f(\x_k,\bar{\y}) \le f(\bar{\x},\bar{\y}) + \epsilon$ for any $k \ge k_2$ and thus 
			\begin{equation}f_{\mu_k}^*(\x_k) \le f(\x_k,\bar{\y})+ \frac{\mu_{k,1}}{2}\|\bar{\y}\|^2 + \mu_{k,2} \le f^*(\x) + 3\epsilon
		\end{equation}for all $k \ge k_2$. By letting $k \rightarrow \infty$, we obtain
		\begin{equation}
		\limsup_{k \rightarrow \infty} f_{\mu_k}^*(\x_k) \le f^*(\bar{\x}) + 3 \epsilon,
		\end{equation}
		and taking $\epsilon \rightarrow 0$ in the above inequality yields the conclusion.
	\end{proof}
	\subsection{Proof of Lemma \ref{lem3}}
	\begin{proof}
		We assume to arrive a contradiction that there exists $\bar{\x} \in \mathbb{R}^m$ satisfying $\x_k \to\bar{\x}$ as $k \to \infty$ with following inequality
		\begin{equation*}
		\liminf_{k \rightarrow \infty} \psi_{\mu_k}(\x_k) < \varphi(\bar{\x}).
		\end{equation*}
		Then, there exist $\epsilon > 0$ and sequences $\x_k \rightarrow \bar{\x} $ and $\{\y_k\}$ satisfying
		\begin{eqnarray}
		&-1 \le f(\x_k, \y_k) - f_{\mu_k}^*(\x_k) \le 0 \quad \label{lem3_eq1} \\
		&\quad F(\x_k, \y_k) \le \psi_{\mu_k}(\x_k) + \epsilon < \varphi(\bar{\x}) - \epsilon.	 \label{lem3_eq1.5}
		\end{eqnarray}
		Then it follows from Lemma \ref{lem2} that 
		\begin{equation}\label{lem3_eq2}
		\limsup_{k \rightarrow \infty} f(\x_k, \y_k) \le \limsup_{k \rightarrow \infty} f_{\mu_k}^*(\x_k) \le f^*(\bar{\x}).
		\end{equation}
		Since either $F(\x,\y)$ or $f(\x,\y)$ is level-bounded in $\y$ locally uniformly in $\bar{\x}$, we have that $\{ \y_k \}$ is bounded. Take a subsequence $\{\y_t\}$ of $\{\y_k\}$ such that $\y_t \rightarrow \hat{\y}$. Then, it follows from Eq. \eqref{lem3_eq1}, Eq. \eqref{lem3_eq2} and the continuity of $f(\x,\y)$ that
		\begin{equation}
		f(\bar{\x},\hat{\y}) \le \limsup_{t \rightarrow \infty} f(\x_t, \y_t) \le f^*(\bar{\x}), 
		\end{equation}
		and thus
		\begin{equation}
		\hat{\y} \in \S(\bar{\x}).
		\end{equation}
		Then, Eq. \eqref{lem3_eq1.5} yields that
		\begin{equation*}
		\varphi(\bar{\x}) \le F(\bar{\x},\hat{\y}) \le \limsup_{k \rightarrow \infty} F(\x_k, \y_k) \le \varphi(\bar{\x}) - \epsilon,
		\end{equation*}
		which implies a contradiction. Thus we get the conclusion.
	\end{proof}
	
	\subsection{Proof of Lemma \ref{lem4}}
	\begin{proof}
		Given any $\x \in \X$, for any $\epsilon > 0$, there exists $\bar{\y} \in \mathbb{R}^n$ satisfying $f(\x,\bar{\y}) \le f^*(\x)$ and $F(\x,\bar{\y}) \le \varphi(\x) + \epsilon$. As $f^*(\x) + \mu_{k,2} \le f_{\mu_k}^*(\x)$, and by the definition of $\varphi_k$, we have
		\begin{equation}
		\begin{aligned}
		\varphi_k(\x) &\le F(\x,\bar{\y}) + \frac{\theta_k}{2}\|\bar{\y}\|^2 - \tau_k\ln\left(f_{\mu_k}^*(\x) - f(\x,\bar{\y})\right) \\
		&\le \varphi(\x) + \epsilon + \frac{\theta_k}{2}\|\bar{\y}\|^2 - \tau_k \ln\mu_{k,2}.
		\end{aligned}
		\end{equation}
		
		By taking $k \rightarrow \infty$ in above inequality, as $\theta_k\rightarrow 0$ and $\tau_k \ln\mu_{k,2} \rightarrow 0$, we have
		\begin{equation}
		\limsup_{k \rightarrow \infty}\varphi_k(\x) \le \varphi(\x) + \epsilon.
		\end{equation}
		Then, we get the conclusion by letting $\epsilon \rightarrow 0$.
	\end{proof}
	
	\subsection{Proof of Proposition \ref{prop2}}
	\begin{proof}
		To prove the epiconvergence of $\varphi_k$ to $\varphi$, we just need to verify that sequence $\{\varphi_k\}$ satisfies the two conditions given in Definition \ref{def_epic}. Considering any sequence $\{\x_k\}$ converging to $\bar{\x}$, since 
		\begin{equation}
		F(\x,\y) \le F(\x,\y) + \frac{\theta_k}{2}\|\y\|^2 - \tau_k\ln\left(f_{\mu_k}^*(\x) - f(\x,\y)\right),
		\end{equation}
		for any $(\x,\y) \in \mathbb{R}^m \times \mathbb{R}^n$ satisfying $-1 \le f(\x,\y)-f_{\mu}^*(\x) < 0$, then we have 
		\begin{equation}
		\psi_{\mu_k}(\x_k) \le \varphi_k(\x_k), \quad \forall k.
		\end{equation}
		By taking $k \rightarrow \infty$ in above inequality, we obtain from Lemma \ref{lem3} that when $\bar{\x} \in \X$,
		\begin{equation}
		\begin{aligned}
		\varphi(\bar{\x}) + \delta_\X(\bar{\x}) &= \varphi(\bar{\x}) \\&\le \liminf_{k \rightarrow \infty} \psi_{\mu_k}(\x_k)\\&
		\le  \liminf_{k \rightarrow \infty} \varphi_k(\x_k) \\ &\le \liminf_{k \rightarrow \infty} \varphi_k(\x_k) + \delta_\X(\x_k).
		\end{aligned}
		\end{equation}
		We have $\liminf_{k \rightarrow \infty} \varphi_k(\x_k) + \delta_\X(\x_k) = + \infty$ when $\bar{\x} \notin \X$, as $\X$ is closed. 
		And thus condition 1 in Definition \ref{epicon} is satisfied. Next, for any $\x \in \mathbb{R}^m$, if $\x \in \X$, then it follows from Lemma \ref{lem4} that
		\begin{equation}
		\limsup_{k \rightarrow \infty} \varphi_k(\x) + \delta_\X(\x) \le \varphi(\x) + \delta_\X(\x).
		\end{equation}
		When $\x \notin \X$, we have $\varphi(\x) + \delta_\X(\x) = + \infty$. 
		Thus condition 2 in Definition \ref{def_epic} is satisfied. Therefore, we get the conclusion immediately from Definition \ref{def_epic}.
	\end{proof}
	
	\subsection{Proof of Theorem \ref{the1}}
	\begin{proof}
		According to Proposition \ref{prop2}, we know that 
		\begin{equation}
		\varphi_k(\x) + \delta_\X(\x) \xrightarrow{e} \varphi(\x) + \delta_\X(\x).
		\end{equation}
		Then the conclusion follows from \cite{AlexanderShapiro2011Perturbation}[Proposition 4.6].
	\end{proof}
	\section{Details of Experiments }\label{sectionC}
	We use PyTorch 1.6 as our computational framework and base our implementation on \cite{grefenstette2019generalized,grazzi2020iteration}.    In all the experiments, we use the Adam method for accelerating the gradient descent of $\x$.  We conducted these experiments on a
	PC with Intel Core i7-9700F CPU, 32GB RAM
	and an NVIDIA RTX 2060S 8GB GPU.                                                                     
	
	\subsection{Toy Examples}
	In numerical experiment, we set $T=100$ for explicit method RHG, $T=100$, $J=20$ for implicit method CG, and $\mu_2=f(x,y)+1$, $(\mu_{k,1}, \theta_k, \tau_k)=(1.0, 1.0, 1.0)/1.01^k $, step sizes $s_1$, $s_2$ and $\alpha$ all equal to $0.01$, $T_{z}=50$, $T_{y}=25$, and $L=1$ in BVFIM.

	\subsection{Hyperparameter Optimization}
	
	In hyperparameter optimization, we set $T=100$ for explicit method RHG, TRHG and BDA, $T=100$, $J=20$ for implicit method CG and Neumann, and $\mu_2=f(\x,\y)$, $(\mu_{k,1}, \theta_k, \tau_k)=(1.0, 1.0, 1.0)/1.01^k $, step sizes $s_1$, $s_2$ and $\alpha$ all equal to $0.01$, $T_{\z}=50$, $T_{\y}=25$, and $L=1$ for BVFIM. We let TRHG truncate at $T/2$ and $\alpha_k=0.5\times 0.999^k$ in BDA.

	We set the training set, validation set, and test set as class balanced.  For each contaminated training sample, we randomly replace its label with a label different from the original one with equal probability. In the calculation of F1 score, if $\x_i\leq0$, we marks the sample $u_i$ as contaminated. In the CIFAR10 experiment, we used an early stop strategy to avoid over-fitting and report the best results achieved. Since $T_{\y}$ gradient descent require $\frac{\partial F}{\partial\y}$ and $\frac{\partial f}{\partial\y}$ separately, we set $T_{\z}+2\times T_{\y}=T$ to fairly compare the time consumed by the algorithm. 

\section{Additional Results}\label{sectionD}

	\subsection{Toy Examples}
 
\begin{figure}  
		\centering  
		\includegraphics[height=3cm,width=4cm]{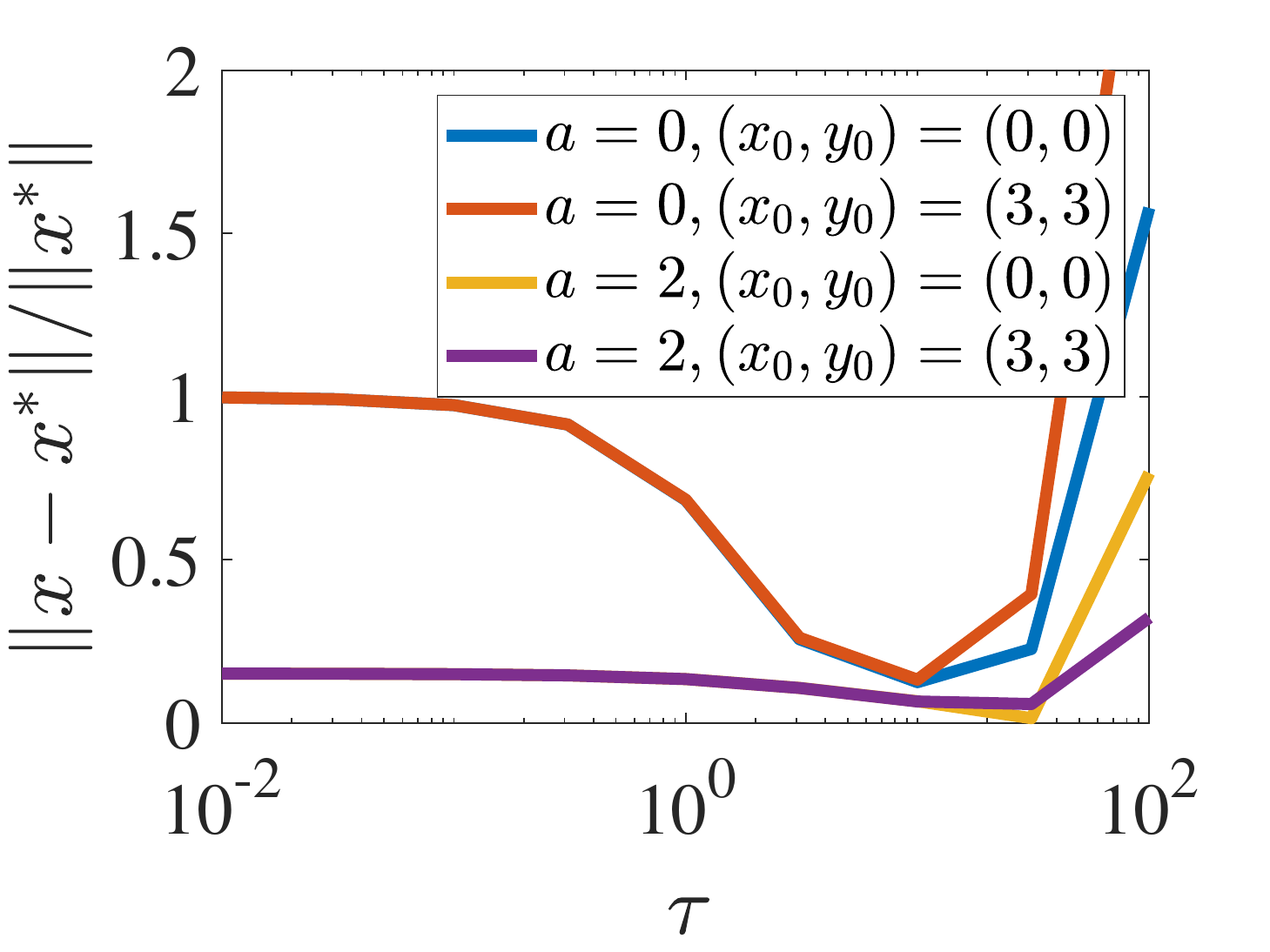}  
		\includegraphics[height=3cm,width=4cm]{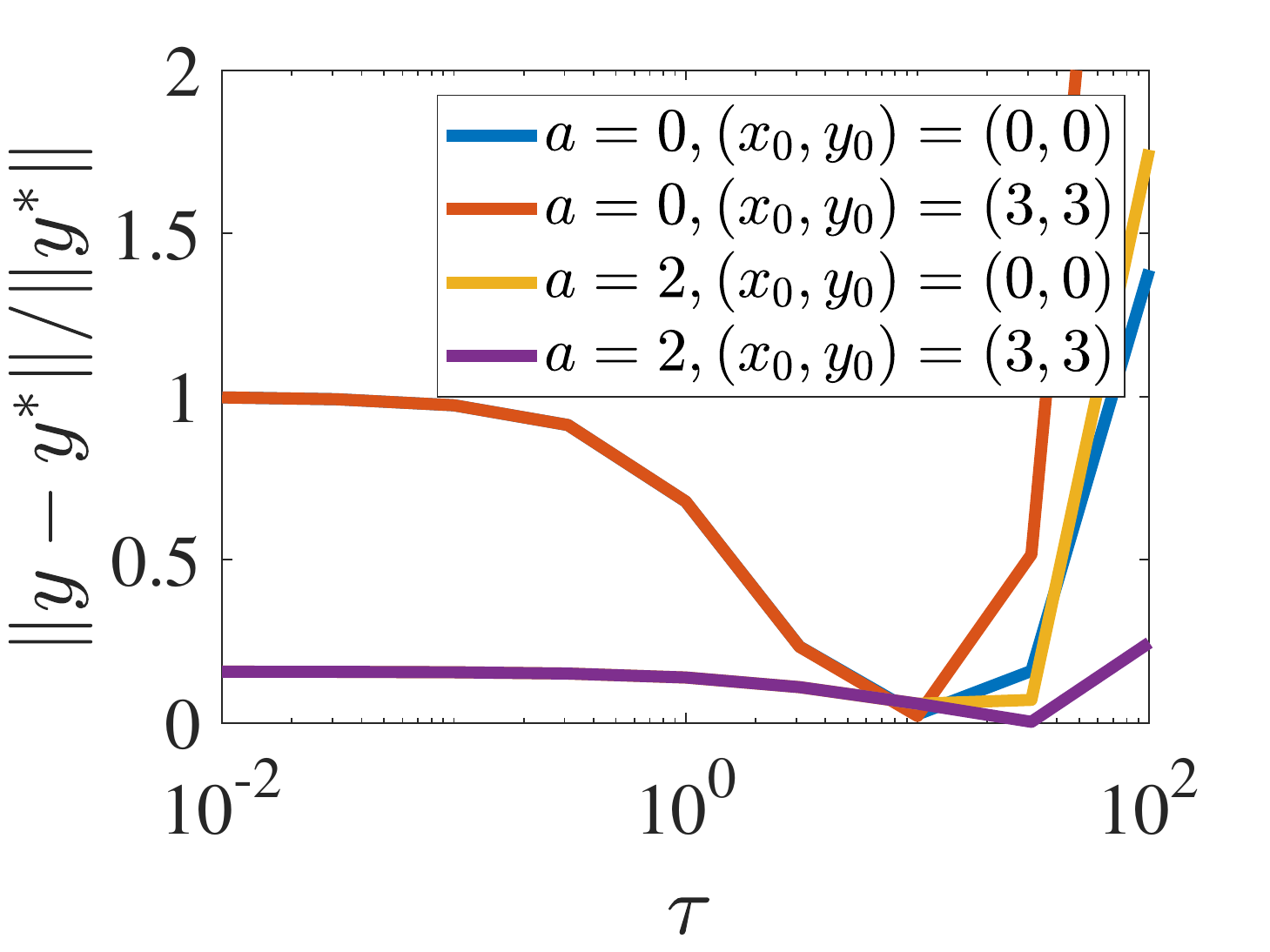}  
		\caption{Convergence results for different regularization coefficients in different initialization settings.}  
		\label{a1}  
	\end{figure}

		We can see that our method has a weaker convergence in LL problem than the existing method under proper initialization (i.e., the initial point is within a locally convex neighborhood of the global optimal point). This is because the main purpose of our experiment is to compare the convergence behavior between different methods and scenarios, so we have not carefully adjusted the regularization coefficients of our methods in order to better show the differences. To verify that our method can also converge as well as the existing methods under proper initialization, we show how to obtain better convergence performance by adjusting $\tau$ in Figure \ref{a1}. It can be seen that an appropriate $\tau$ can greatly improve the convergence behavior. In addition, we validate this with a larger LLC problem in the section \ref{B3}.
	\begin{figure}  
		\centering  
		\includegraphics[height=3cm,width=4cm]{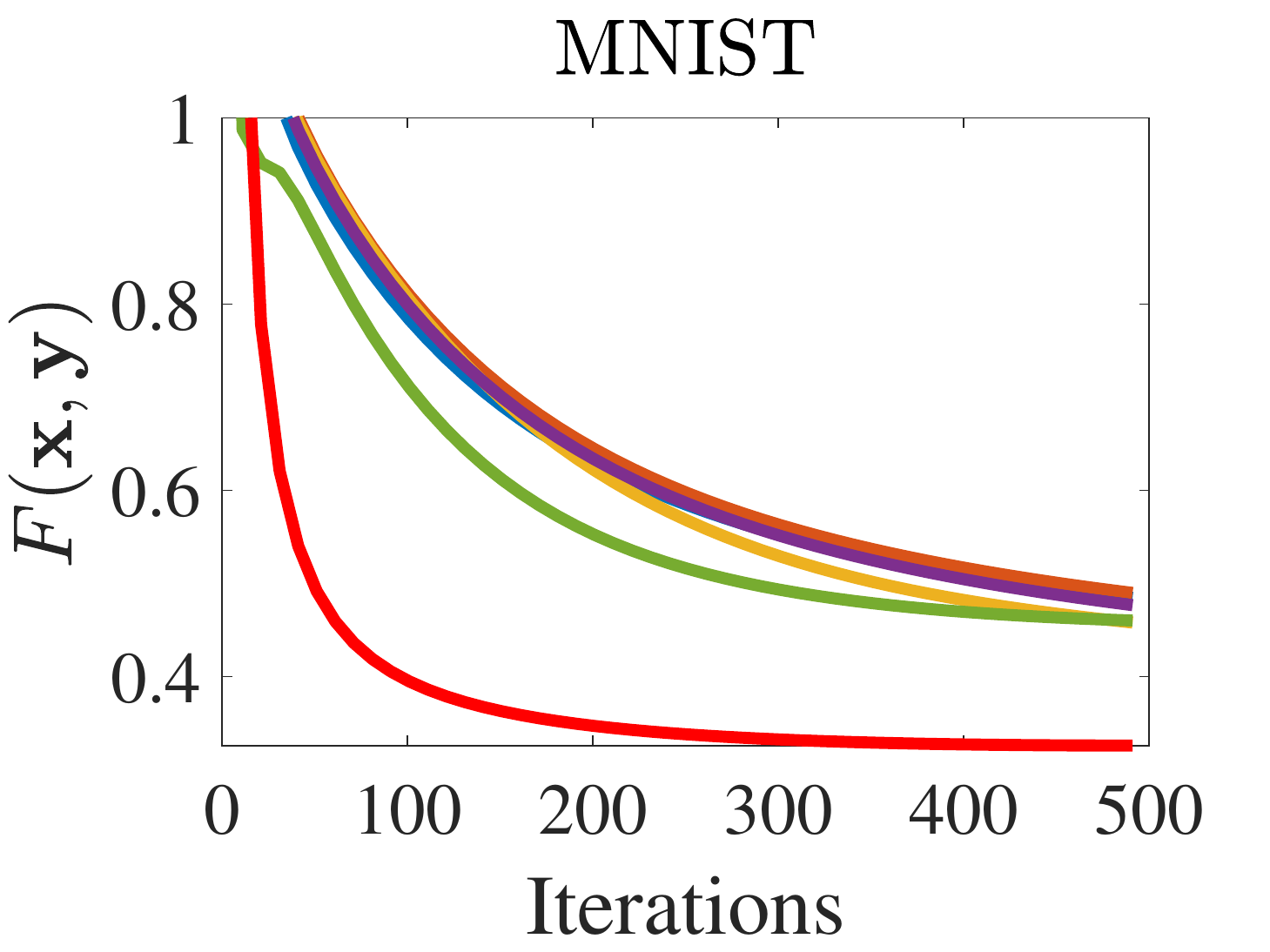}  
		\includegraphics[height=3cm,width=4cm]{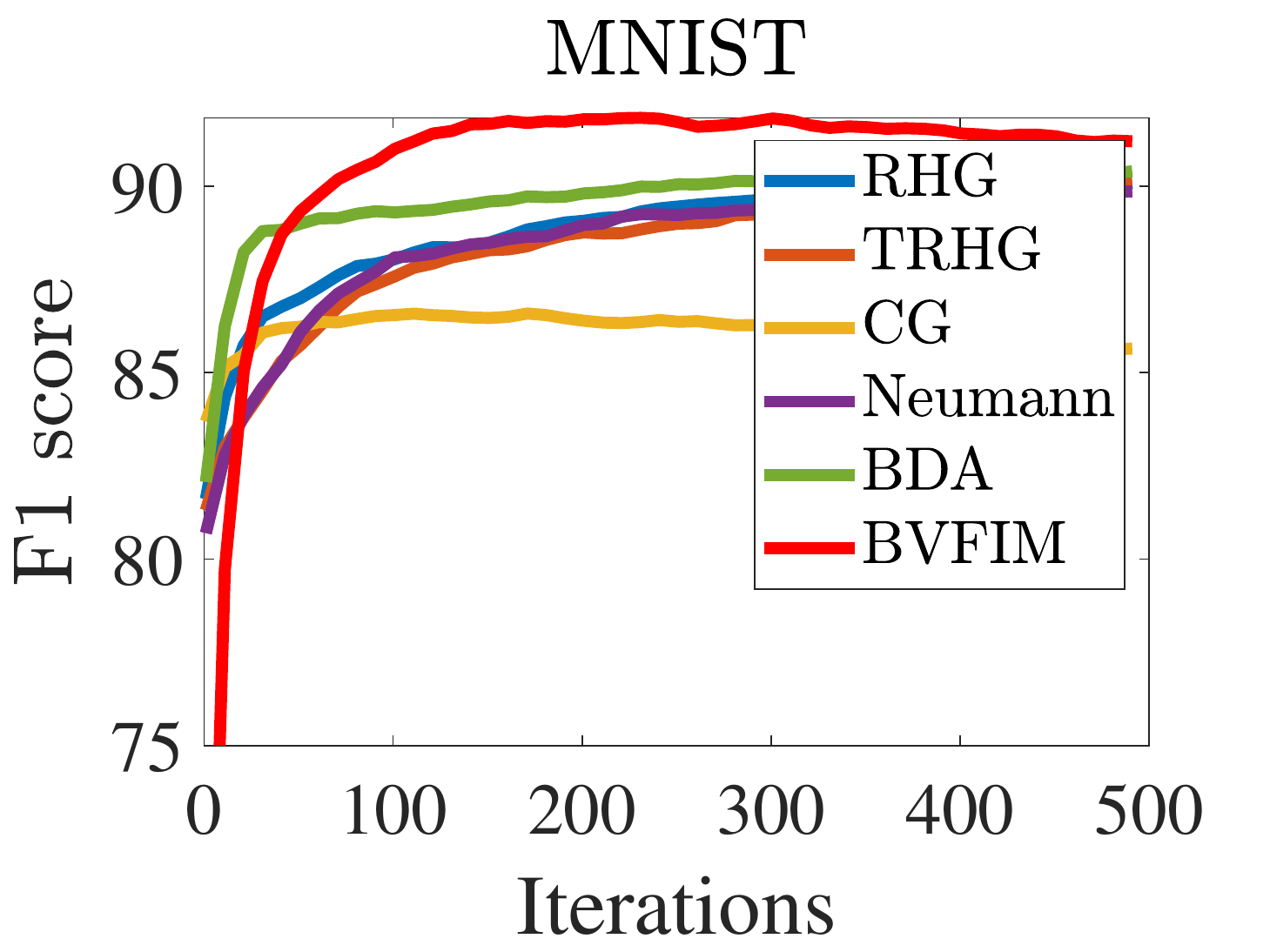}  
		\includegraphics[height=3cm,width=4cm]{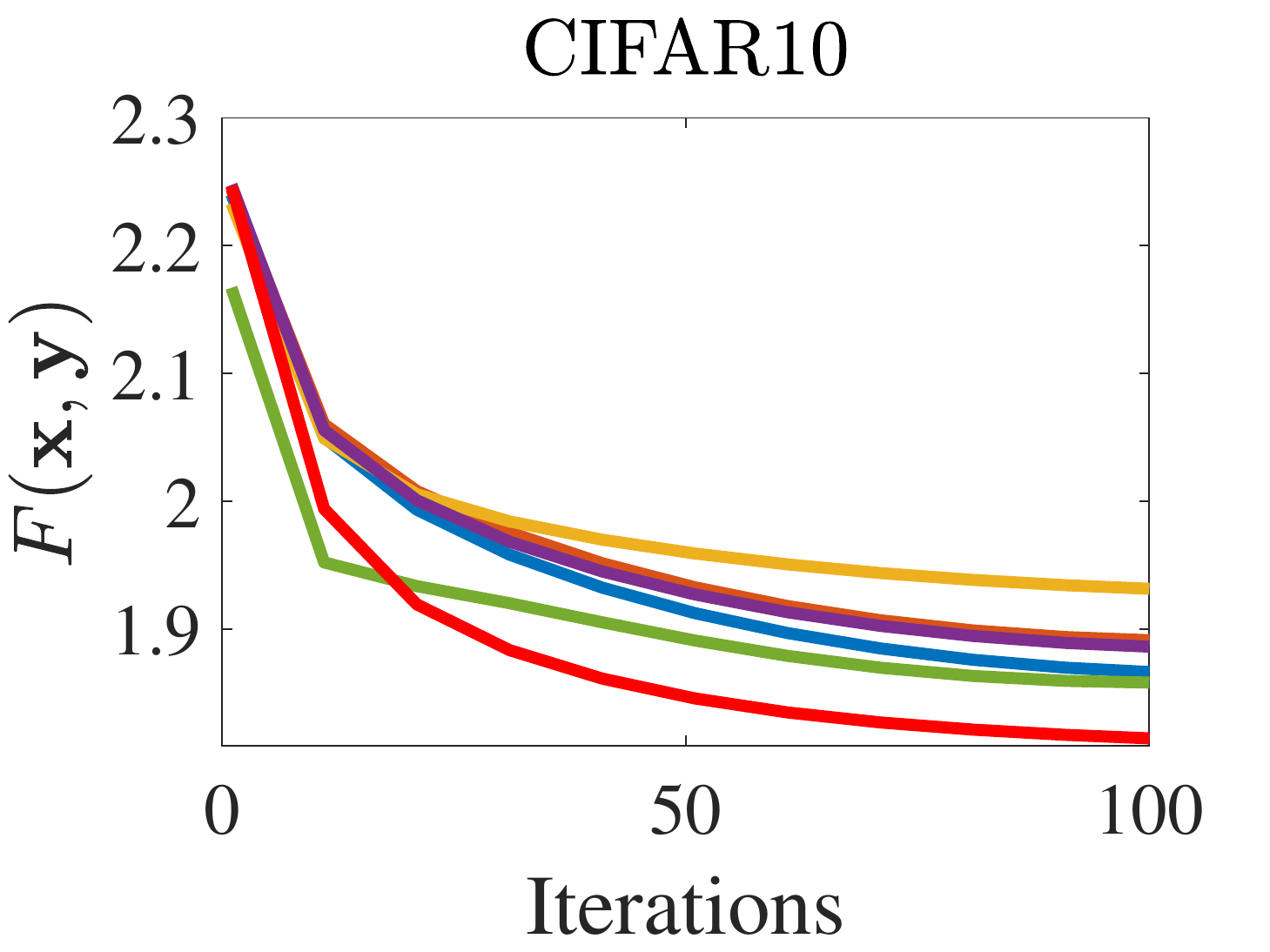}  
		\includegraphics[height=3cm,width=4cm]{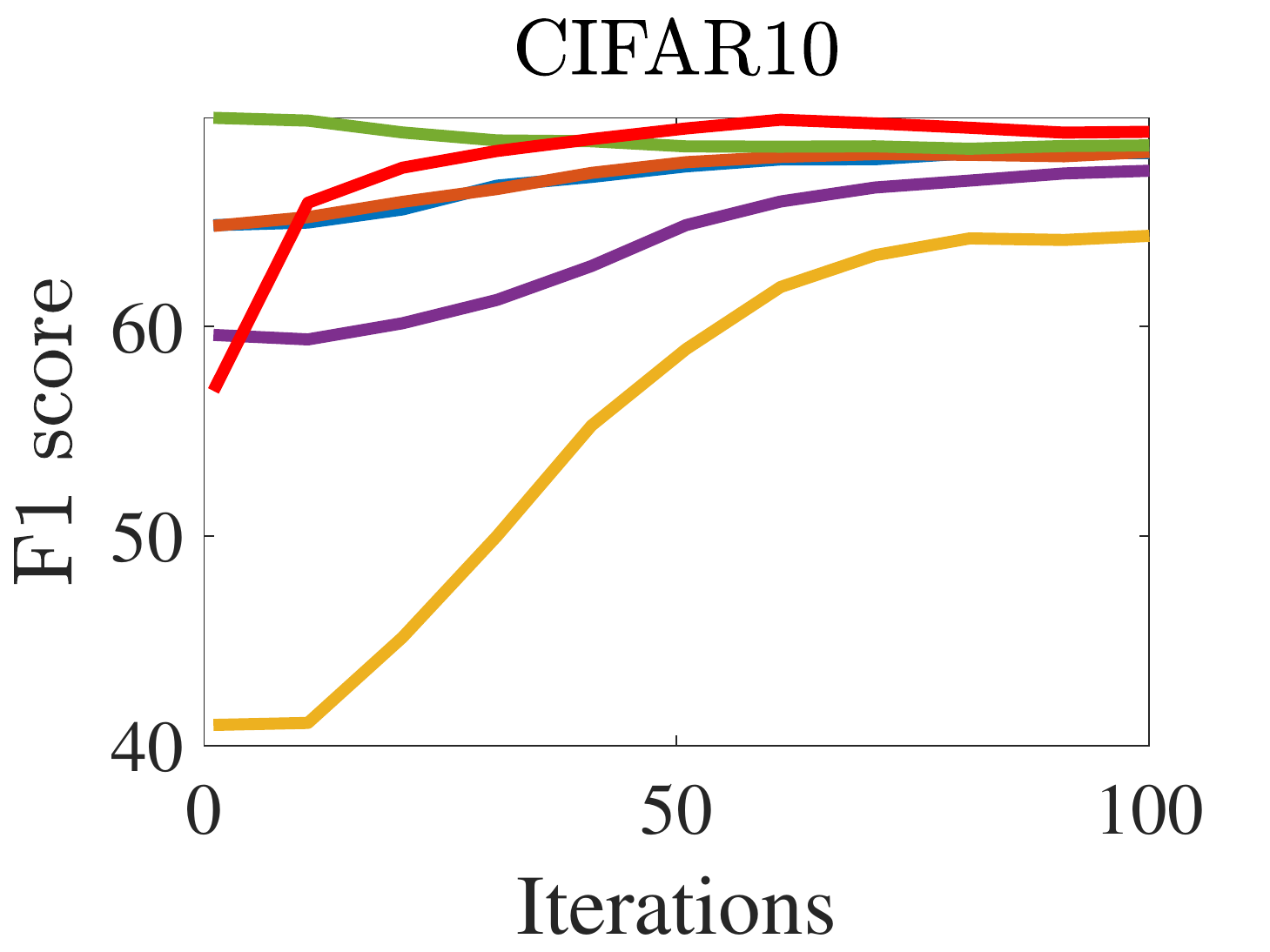}  
		\caption{$F(\x,\y)$ and F1 score between existing methods and BVFIM. The curves are based on the MNIST and CIFAR10 experiment. The legend is only plotted in the second subfigure.}  
		\label{a2}  
	\end{figure} 
	\subsection{Hyperparameter Optimization}
	The UL objective and F1 scores of BVFIM and compared methods on the MNIST and CIFAR10 dataset are plotted in Figure \ref{a2}.

\begin{figure}
			\centering  
			\includegraphics[height=3cm,width=4cm]{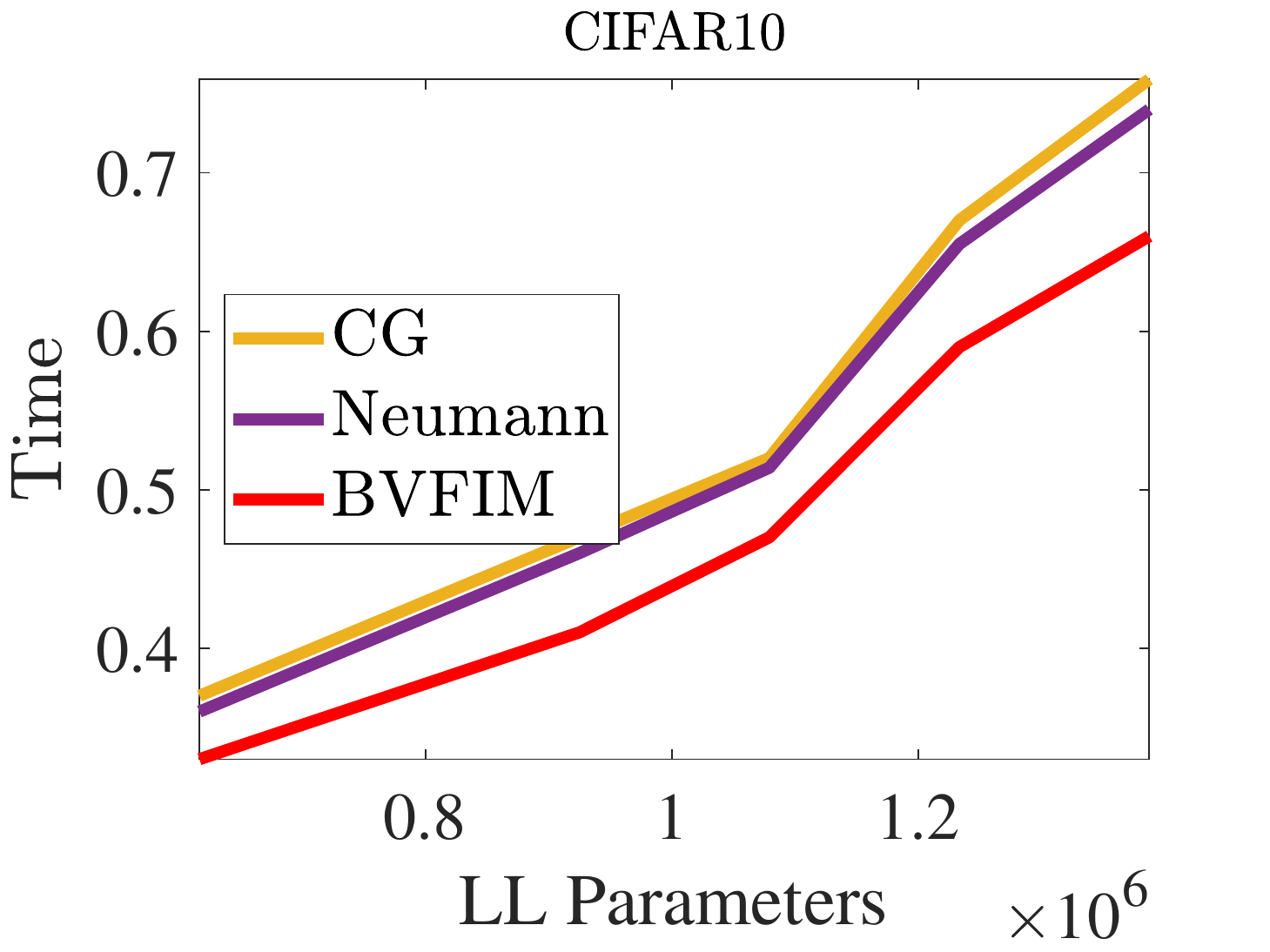}  
			\includegraphics[height=3cm,width=4cm]{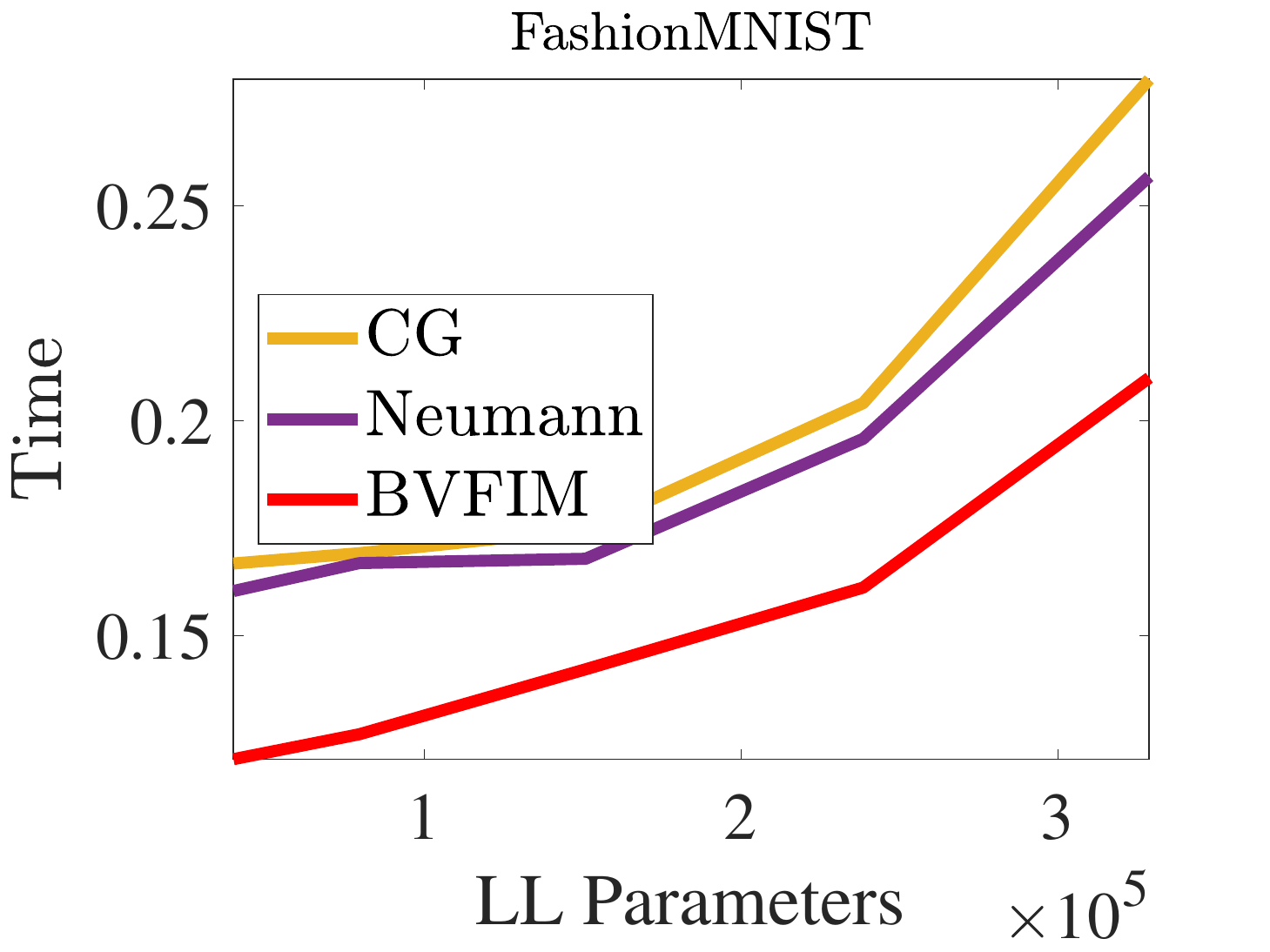} 
			\caption{Comparison of calculation time of IGBMs and BVFIM under different LL parameter quantities.}  
			\label{net}
		\end{figure}

	In addition, we verify the computation time variation of BVFIM and existing gradient-based methods under different LL variable dimensions  on the FashionMNIST and CIFAR10 dataset. In order to show the comparison results more clearly, we compared with IGBMs which are faster in the existing gradient-based methods.
	Figure~\ref{net} shows the computation time with different LL variable parameter quantities. It can be seen that BVFIM is faster than IGBMs at different parameter quantities, and this advantage becomes more significant as the number of LL parameters increases.
	\subsection{Additional LLC Experiments}\label{B3}

	To verify the validity of the BVFIM method in conventional LLC problems, we supplemented a meta-learning experiment.
	The goal of meta-learning is to learn an algorithm that can handle new tasks well. In particular, we consider the few-shot learning problem, where each task is a N-way classification and it aim to learn the hyperparameter $\x$ so that each task can be solved by only M training samples. (i.e. N-way M-shot)

	Similar to recent work~\cite{franceschi2018bilevel,liu2020generic}, we modeled the network in two parts: a four-layer convolution network $\x$ as a common feature extraction layer between different tasks, and logical regression layer $\y={\y^i}$ as separate classifier for each task. We also set dataset as $\D=\{\D^i\}$, where $\D^i=\D^i_{\mathtt{tr}}\cup \D^i_{\mathtt{val}}$ is for the $i$-th task. Then we set the loss function of the $j$-th task to $\mathtt{CE}(\x,\y^i;\D^i_{\mathtt{tr}})$ for the LL problem, thus the
	LL objective can be defined as
	\begin{equation*}
	f(\x,\y)=\sum_i \mathtt{CE}(\x,\y^i;D^i_{\mathtt{tr}})
	\end{equation*}
	As for the UL objective, we also utilize cross-entropy function but define it based on $\{D^i_{\mathtt{val}}\}$ as
	\begin{equation*}
	F(\x,\y)=\sum_i \mathtt{CE}(\x,\y^i;D^i_{\mathtt{val}})
	\end{equation*}
	Our experiment was performed on two widely used benchmark datasets: Omniglot~\cite{lake2015human}, which contains examples of 1623 different handwritten characters from 50 alphabets and MiniImagenet~\cite{vinyals2016matching}, which is a subset of ImageNet~\cite{deng2009imagenet} that contains 60000 downsampled images from 100 different classes. We compared our BVFIM to several approaches, such as RHG, TRHG and BDA~\cite{liu2020generic}.
	\begin{table}
	\caption{The averaged few-shot classification accuracy on Omniglot and MiniImageNet (M=1)}
	\begin{center}
		\begin{small}
			\begin{tabular}{lccc}
				\toprule\label{T21}
				\multirow{2}{*}{Method}&\multicolumn{2}{c}{Omniglot}&MiniImagenet\\
				&5-way&20-way&5-way\\
				\midrule
				RHG&98.60&95.50&48.89\\
				TRHG&98.74&95.82&47.67\\
				BDA&\textbf{99.04}&\textbf{96.50}&49.08\\
				BVFIM&98.85&95.55&\textbf{49.28}\\
				\bottomrule
				
				\end{tabular}
			\end{small}
		\end{center}
	\end{table}
	
	For RHG, TRHG and BDA, we follow the settings in ~\cite{liu2020generic}. For BVFIM, we set $T_{\z}=5$, $T_{\y}=10$, $\mu_2=f(\x,\y)$, $(\mu_{k,1}, \theta_k, \tau_k)=(1, 0.1, 10)/k $, step sizes $(s_1,s_2,\alpha)=(0.01,0.01,0.001)$, $L=1$ and $K=50000$. We set meta-batch size of 16 episodes for Omniglot dataset and of 4 episodes for Miniimagenet dataset. We set $\y_{k,l}^0=\z_{k,l}^{T_{\z}}$ to warm up. 
	
	It can be seen in Table \ref{T21} that BVFIM can get slightly poorer performance than existing methods on Omniglot dataset and get the best performance on MiniimageNet dataset, which proves that our method can also obtain competitive results for LLC problems.

\end{document}